\documentclass{article}
\usepackage{graphicx} 
\usepackage{amsmath,amsfonts,amssymb,amsthm}
\usepackage{xcolor}
\usepackage{hyperref}

\newtheorem{remark}{Remark}[section]

\newcommand{\R}{\mathbb{R}}

\title{Control of Overpopulated Tails in  Kinetic Epidemic Models}

\date{\today}

\author{M. Zanella \thanks{Department of Mathematics "F. Casorati", University of Pavia, Italy. Email:mattia.zanella@unipv.it} \and A. Medaglia\thanks{Mathematical Institute, University of Oxford, United Kingdom. Email:andrea.medaglia@maths.ox.ac.uk}}
\begin{document}

\maketitle

\begin{abstract}
We introduce model-based transition rates for controlled compartmental models in mathematical epidemiology, with a focus on the effects of control strategies applied to interacting multi-agent systems describing contact formation dynamics. In the framework of kinetic control problems, we compare two prototypical control protocols: one additive control directly influencing the dynamics and another targeting the interaction strength between agents. The emerging controlled macroscopic models are derived for an SIR compartmentalization to illustrate their impact on epidemic progression and contact interaction dynamics.  Numerical results show the effectiveness of this approach in steering the dynamics and controlling epidemic trends, even in scenarios where contact distributions exhibit an overpopulated tail. 
\end{abstract}

\section{Introduction}

In this work, we focus on the derivation of model-based transition rates for controlled compartmental models in mathematical epidemiology, emphasizing the impact of a centralized control strategy applied to interacting multiagent systems describing contact formation dynamics. Control methods in compartmental epidemiology have been extensively studied to define optimal model-oriented non-pharmaceutical interventions over variable time horizons \cite{APZ0,APZ1,Ben,CFCGCC,CG2,LCCC,WAP}. Further approaches have focused on control methods to optimize vaccination procedures and quarantine measures, see, e.g., \cite{AM,BSBG,Brauer,CG,CMR,Ferguson}. In this direction, the problem of controlling  multiagent systems has recently emerged as a natural follow-up to the description and modeling of their self-organizing features \cite{AP,BBCK,CFPT,CKRT,CKPP}. In the case of large systems, significant contributions have been made in the context of the derivation of control strategies for mean-field and kinetic equations, see, e.g., \cite{ACFK,AHP,APTZ,FPR,FS,TZ} and the references therein. 

The evolution of epidemics can be understood as the result of unobservable interactions among a large number of agents of heterogeneous populations, whose transitions are influenced by multidimensional factors, including both behavioural aspects and the biological characteristics of the pathogen \cite{De,Dong,Eubank,Fumanelli,GIV,LT,SKLLFJ}. In this context, kinetic theory emerges as a robust mathematical framework to understand the emerging statistical properties of large multiagent systems \cite{BL,BZ,CPS,DPaTZ,DPTZ,Z2}. A key motivation for introducing kinetic-type equations lies in their multiscale features, which connects the microscopic level of individual agents, where fundamental interactions occur, to the macroscopic level, where data are typically available. This characteristic links the kinetic approach to the effective evolution of multiagent systems. Moreover, the ability to derive macroscopic models that are consistent with the underlying microscopic interaction dynamics is crucial for parameter estimation, which in turn supports effective decision-making processes. In this context, we refer to \cite{DTZ,FMZ,Z_etal} for comparisons of these modelling approaches with structured real data, particularly with respect to demographic factors such as age.

In this work, we focus on a kinetic model for epidemic dynamics in which the population of agents is divided into compartments and the mass of each compartment changes in connection with microscopic dynamics defining the emerging contact distribution of the system which characterises the heterogeneity of the population \cite{B_etal,BL,MA,Zha}. The dynamics of the particles are governed by a simple transition operator consistent with growth-type operators  \cite{PTZ,Z2}. At the kinetic level, this model may lead to the emergence of an equilibrium distribution characterized by power-law tails. This behaviour is associated with an increased probability of encountering agents with a large number of contacts, and consequently, with a higher potential for epidemic spread. To highlight this phenomenon we considered a macroscopic model derived in the presence of a local incidence rate which depends, at variance with the approach defined in \cite{DPTZ,DTZ}, on higher-order moments of the contact distribution.

To tackle this scenario, we propose a control mechanism acting on the contact dynamics, with the goal of converting a fat-tailed distribution into a slim-tailed distribution.  Existing strategies, although effective in steering the mean number of contacts toward a target value, appear insufficient to contain epidemic transmission in the presence of overpopulated tails.

{In more detail, we consider two control strategies aimed at influencing the dynamics of the multiagent system. The first approach employs an additive control mechanism that directly modifies the agents' individual dynamics, while the second targets the interaction strength between agents, effectively shaping the overall interaction structure. By deriving the corresponding controlled macroscopic models, we apply these strategies to an SIR compartmentalization of the society. This allows us to investigate their impact on both the progression of an epidemic and the dynamics of contact interactions within the population. The comparison highlights that controlling the interaction is capable of influencing the nature of the contact structure of agents by shaping the tail of the contact distribution. The effectiveness of the control is tested by means of the direct computation of the emerging equilibria, whose moments evolve through a consistent macroscopic model obtained through an  equilibrium closure. }

The present paper is organised as follows. In Section \ref{sect:2} we recall basic concepts of kinetic modelling for epidemic dynamics, and the emerging equilibrium densities of contact formation dynamics have been obtained. Hence, in Section \ref{sect:22}-\ref{sect:23} we embed the obtained equilibrium density in the kinetic model to derive the evolution of observable quantities. The control of the tail behaviour is discussed in Section \ref{sect:3}  by comparing two control strategies. In Section \ref{sect:31} the emerging distributions of the resulting Fokker-Planck equations for the two control protocols are computed and compared, and the macroscopic equations are obtained in Section \ref{sect:32}. Finally, in Section \ref{sect:4} we present several numerical results.  

\section{Kinetic compartmental models with contact structure}   \label{sect:2}

We consider a prototypical system of agents subdivided into the following epidemiological relevant states: susceptible ($S$) agents are the ones that can contract the disease, infected and infectious, ($I$) agents are responsible for the spread of the disease, and removed ($R$) agents cannot spread the disease. In the following, we will indicate with $\mathcal C = \{S,I,R\}$ the considered compartmentalization of the population. 

We denote by $f_J = f_J(x,t):\R^+\times \R^+ \to \R^+$, $J \in \mathcal C$, the distribution of the number of contacts $x\in\mathbb R^+$ at time $t\ge0$ of agents in the compartment $J\in \mathcal C$. Hence, the contact distribution $f(x,t)$ of the society is obtained as follows
\[
\sum_{J \in \mathcal C} f_J(x,t) = f(x,t), \qquad \int_{\mathbb R^+}f(x,t) dx = 1,
\]
and the mass fractions of the population in each compartment are defined as
\[
\rho_J(t) = \int_{\mathbb R^+} f_J(x,t)dx, 
\]
while their moment of order $r>0$ are given by
\[
\rho_J(t)m_{r,J}(t) = \int_{\mathbb R^+}x^r f_J(x,t)dx.
\]
In the following, we discuss the evolution of the contact distributions $f_J(x,t)$ due to microscopic interactions.
In the following, to simplify notations we will indicate with $m_J(t) = m_{1,J}(t)$ the mean, corresponding then to the cases $r=1$.

\subsection{Kinetic models for epidemic spreading}\label{sect:22}
Following the approach presented in \cite{DPTZ}, we are interested in the evolution of the kinetic densities $(f_J)_{J \in \mathcal C}$ solution to 
\begin{equation}
\label{eq:kinetic}
    \begin{cases}
\partial_t f_S(x,t) = -K(f_S,f_I)(x,t) + \dfrac{1}{\tau}Q_S(f_S)(x,t), \\
\partial_t f_I(x,t) =K(f_S,f_I)(x,t)-\gamma_I f_I(x,t) + \dfrac{1}{\tau}Q_I(f_I)(x,t) \\
\partial_t f_R(x,t) = \gamma_I f_I(x,t) +  \dfrac{1}{\tau}Q_R(f_R)(x,t),
    \end{cases}
\end{equation}
where $\gamma_I>0$ is the recovery rate, $\tau>0$ is the time scale of the contact dynamics, and the $Q_J$ are suitable thermalisation operators that will be introduced in Section \ref{sect:contact}. In \eqref{eq:kinetic}, the transmission of the infection is governed by the local incidence rate
\begin{equation}
\label{eq:incidence}
K(f_S,f_I)(x,t) = f_S(x,t) \int_{\mathbb R^+} \kappa(x,x_*) f_I(x_*,t)dx_*,
\end{equation}
being $\kappa(x,x_*)$ the contact function weighting the frequency of contacts between susceptible and infected agents. We highlight that the case $\kappa(x,x_*) \equiv \beta>0$ corresponds to the hypothesis of a homogeneous transmission between agents as in the classical SIR model. To take into account the influence of the shape of the contact distribution of infected in the transmission dynamics we may consider
\begin{equation}
\label{eq:kappa}
\kappa(x,x_*) =  \sum_{\ell=1}^L \beta_\ell (x x_*)^{\ell}, 
\end{equation}
where $\beta_\ell\ge0$ are proportionality parameters, and $\ell=1,\dots,L$ determines the relevant moments of the contact distribution in the infection dynamics. In \cite{DPTZ}, the case $L = 1$ has been considered, for which the mean number of connections of the infected population weights the spread of the infection. In this case we get
\[
K(f_S,f_I)(x,t) = \beta_1 x f_S(x,t) \rho_I(t) m_{I}(t).
\]
More generally, higher order moments of the distribution are necessary to characterize the spread of the infection. To understand the impact of the tails of the contact distribution, we can  consider the case $L=2$. Hence, under this hypothesis the incidence rate defined in \eqref{eq:incidence} reads
\[
K(f_S,f_I)(x,t) = \beta_1 x f_S(x,t) \rho_I(t) m_{I}(t) + \beta_2 x^{2}f_S(x,t) \rho_I(t) m_{2,I}(t). 
\]
In this case, we obtain information on the impact of the mean and energy of the contact distribution on the evolution of the disease. 
Within the choice in \eqref{eq:kappa}, the incidence rate is proportional to the product of the mean number of contact of susceptible and infected agents and on their energies. 

\subsection{Contact formation dynamics}\label{sect:contact}

We consider a large system of $N\gg0$ agents that are identified by their number of daily contacts $x_i \in \R^+$, $i = 1,\dots,N$. The dynamics of each agent is modified through deterministic variations of the contact structure due to environmental factors, and by random fluctuations in the contact variation due to population heterogeneities. Following the approach presented in \cite{DPTZ} and indicating by $x^\prime$ the post-interaction number of contacts of an agent, we get
\begin{equation}
\label{eq:micro}
x^\prime = x + \epsilon \Psi(x/m_J) x + x \eta,
\end{equation}
where 
\begin{equation} \label{eq:psi}
\Psi(x/m_J) = \frac{\alpha}{2\delta} \left( \left( \frac{x}{m_J}\right)^\delta - 1\right)
\end{equation}
is a generic growth function which describes logistic-type growths for any $\delta>0$, von Bertalanffy-type growths for $\delta<0$ and coincides with a  Gompertz growth function in the limit $\delta \to 0$, see \cite{PTZ}. In \eqref{eq:micro}, we further introduced the random variable $\eta$ such that, denoting by $\left\langle \cdot \right\rangle$ the expectation with respect to the distribution of the introduced random variable, we get $\left\langle \eta \right\rangle = 0$ and $\left\langle \eta^2 \right\rangle = \epsilon \sigma^2$, while in \eqref{eq:psi} $\alpha>0$ is a proportionality constant.

Hence, following the approach in \cite{Cerc,PT}, we can characterise the contact formation dynamics for the generic compartment $f_J=f_J(x,t)$ by considering the following space homogeneous Boltzmann-type equation 
\begin{equation} \label{eq:boltz}
\frac{\partial }{\partial t} f_J(x,t) = Q_J(f_J)(x,t) := \left\langle\int_{\R^+} B(x)  \dfrac{1}{{}^\prime \mathcal{J}} f_J({}^\prime x,t) dx \right\rangle - f_J(x,t),
\end{equation}
where ${}^\prime x\ge0$ denote the pre-transition number of contacts that produces the post-transition contacts $x\ge0$ following the transition scheme \eqref{eq:micro}, while ${}^\prime \mathcal{J}$ denotes the determinant of the Jacobian of the transformation ${}^\prime x \to x$. In \eqref{eq:boltz}, the term $B(x):\R^+\to\R^+$ is the collision kernel tuning the interaction frequency which may depend on the number of contacts $x\ge0$. In \cite{FMZ}, the kernel
\begin{equation}\label{eq:kernel}
B(x) = x^{-\frac{1+\delta}{2}},
\end{equation}
has been considered. 
The kinetic equation \eqref{eq:boltz} can be conveniently written in weak form by introducing a test function $\varphi(\cdot)$. Hence, the evolution of observable quantities is given as follows
\[
\dfrac{d}{dt} \int_{\R^+}\varphi(x)f_J(x,t)dx = \int_{\R^+} B(x)\left \langle \varphi(x^\prime) - \varphi(x) \right \rangle f_J(x,t) dx.
\]
The introduced operator is always mass preserving: by choosing $\varphi(x) = 1$ we get
\[
\int_{\mathbb R^+} Q_J(f_J)(x,t) dx = 0.
\]
It is momentum-preserving if $\delta \equiv \pm 1$ since from \eqref{eq:kinetic} we get
\[
\int_{\mathbb R^+} x Q_J(f_J)(x,t) dx = -\dfrac{\alpha}{2\delta} \int_{\mathbb R^+} x^{1-\frac{1+\delta}{2}} \left( \left( \dfrac{x}{m_J}\right)^{\delta}-1 \right)f_J(x,t)dx. 
\]
More generally, for any $\delta \in (-1,1)$ the momentum is not a conserved quantity for the introduced collision operator.

In view of the complexity of obtaining an analytical expression on the large time behaviour of the kinetic equation \eqref{eq:boltz}, it is possible to derive from the introduced kinetic model a reduced complexity Fokker-Planck-type partial differential equation for which the study of asymptotic properties is easier, see \cite{PT}. In the quasi-invariant limit, we get the following operator 
\begin{equation}
\label{eq:FPQJ}
\bar{Q}_J(f_J)(x,t) = \dfrac{\alpha}{2\delta}\partial_x \left[ x^{1-\frac{1+\delta}{2}} \left( \left(\dfrac{x}{m_J} \right)^\delta-1 \right)f_J + \dfrac{\sigma^2}{2}\partial_x \left(x^{2-\frac{1+\delta}{2} }f_J\right) \right]. 
\end{equation}
For the sake of brevity, we omit here the details of the derivation, we point the interested reader to \cite{DPTZ} for more details. We can observe how the derived Fokker-Planck operator $\bar Q_J(\cdot)$ is mass preserving and momentum preserving if $\delta = \pm 1$. 

\begin{remark}
The process defined in \eqref{eq:micro} can be obtained in terms of a system of stochastic differential equations (SDEs) for the dynamics of each agent of the system $\{x_i\}_{i=1}^N$ where the dynamics characterising the number of contacts of each agent is described as follows
\[
dx_i = B(x_i)\Psi(x_i)dt + \sqrt{2\sigma^2 x_i}dW_i^t,
\]
being $\{W_i^t\}_{i=1}^N$ a set of independent Wiener processes and $B(x_i)\ge0$ a factor describing the interaction frequency of each agent, see \cite{Z2}. Following \cite{LeBL} it can be shown that in the limit $N\to +\infty$ the evolution of the system of agents is obtained in terms of the evolution of a probability density solution to the Fokker-Planck equation \eqref{eq:FPQJ}. 
\end{remark}

\subsubsection{Equilibrium distribution of the Fokker-Planck equation }
\label{subsect:equilibrium}
The Fokker-Planck-type operator defined in \eqref{eq:FPQJ} is such that its equilibrium density $f^q_{J}(x)$ can be computed as the unique solution to the following differential equation 
\[
\dfrac{\alpha}{\delta} x^{1-\frac{1+\delta}{2}} \left(\left( \dfrac{x}{m_J} \right)^{\delta}-1 \right)f_J^q(x) + \sigma^2 \partial_x (x^{2-\frac{1+\delta}{2}}f_J^q(x)) = 0, \qquad \delta \in [-1,1]
\]
which is given by
\begin{equation}
\label{eq:fqJ}
f_J^q(x) = C_{\delta,\sigma^2,\alpha} x^{\frac{\alpha}{\sigma^2\delta}-2+\frac{1+\delta}{2}}\exp\left\{-\dfrac{\alpha}{\sigma^2\delta^2}\left( \dfrac{x}{m_J} \right)^{\delta} \right\},
\end{equation}
where $C_{\delta,\sigma^2,\alpha}>0$ is a normalization constant. 


The equilibrium solution obtained in \eqref{eq:fqJ} depends on the parameter $\delta \in [-1,1]$ of the model which influences the behaviour of the tails. Indeed, we may observe how, for any $0\le\delta\le 1$ the equilibrium distribution is characterized by slim tails since it exhibits exponential decay for $x\gg0$. In particular, in the case $\delta = 1$ we get a Gamma distribution
\begin{equation}
\label{eq:gamma}
f^q_J(x) = \dfrac{\lambda^\lambda}{m_J^\lambda \Gamma(\lambda)}x^{\lambda-1}\exp\left\{-\dfrac{\lambda x}{m_J}\right\}, \qquad \lambda = \dfrac{\alpha}{\sigma^2},
\end{equation}
having slim tails. In the case of the Gamma distribution, we have the following relations on the first and second order moment
\[
\int_{\R^+} xf^q_J(x)dx = m_J, \qquad  \int_{\R^+} x^2f^q_J(x)dx = \dfrac{\lambda+1}{\lambda}m_J^2.
\]
Furthermore, in the special case $\delta \to 0$ it can be shown that the equilibrium density is a lognormal distribution \cite{PTZ}. On the other hand, for any $\delta <0$ the equilibrium distribution is a fat tails distribution with polynomial decay for $x\gg 0$. In particular, if $\delta = -1$, we obtain the inverse Gamma distribution 
\begin{equation}
\label{eq:invgamma}
f^q_J(x) = \dfrac{(\lambda m_J)^{\lambda+1}}{\Gamma(\lambda+1)}x^{-2-\lambda}\exp\left\{ -\dfrac{\lambda m_J}{x}\right\}, \qquad \lambda = \dfrac{\alpha}{\sigma^2}, 
\end{equation}
for which  we have
\[
\int_{\R^+} xf^q_J(x)dx = m_J, \qquad  \int_{\R^+} x^2f^q_J(x)dx = \dfrac{\lambda}{\lambda-1}m_J^2. 
\]
Therefore, the sign of the generally unknown parameter $\delta \in [-1,1]$ have a strong impact in terms of the behaviour of the multi-agent system. 



\subsection{Macroscopic equations}  \label{sect:23}

In this section, we derive the evolution of observable quantities of the kinetic system of equations \eqref{eq:kinetic}. In Section \ref{subsect:equilibrium} we have observed how the conserved quantities of the contact formation dynamics are given by mass and momentum, for which we wish to obtain a closed system of equations describing their evolution. In the following, we introduce an approximation which is reminiscent of the
equilibrium closure in classical kinetic theory, see e.g. \cite{BGL}. 

As observed in \cite{Z2}, under the hypothesis $\kappa(x,x_*) \equiv \beta$, corresponding to the case where the infection dynamics is independent from the contact distribution, we obtain a system defining the evolution of mass fractions given by the classical SIR model. Indeed, integrating \eqref{eq:kinetic} in $\R^+$ we get
\[
\begin{split}
\dfrac{d}{dt}\rho_S &= -\beta\rho_S \rho_I,  \\
\dfrac{d}{dt} \rho_I &= \beta \rho_S\rho_I - \gamma_I \rho_I, \\
\dfrac{d}{dt} \rho_R &= \gamma_I \rho_I
\end{split}\]
whereas from \eqref{eq:kinetic} we may compute $1/\rho_J(t)\int_{\R_+}xf_J(x,t)dx$ from which we get $\frac{d}{dt}m_J(t) = 0$ for all $J\in \mathcal C$. A first attempt to incorporate the impact of social contacts in the dynamics is obtained by choosing $L = 1$ from which we get
\begin{equation}
\label{eq:mass_L1}
\begin{split}
\dfrac{d}{dt}\rho_S &= - \beta_1\rho_S m_{S} \rho_I m_{I},  \\
\dfrac{d}{dt} \rho_I &= \beta_1\rho_S m_{S}\rho_I m_{I} - \gamma_I \rho_I, \\
\dfrac{d}{dt} \rho_R &= \gamma_I \rho_I,
\end{split}
\end{equation}
such equation can be coupled with the evolution of the mean number of connections which is a conserved quantity of the collision operator \eqref{eq:FPQJ}. We get 
\begin{equation}
\label{eq:mom_L1}
\begin{split}
\dfrac{d}{dt}\rho_S m_{S} &= -\beta_1\rho_S m_{2,S} \rho_I m_{I},  \\
\dfrac{d}{dt} \rho_I m_{I} &= \beta_1 \rho_S m_{2,S}\rho_I m_{I} - \gamma_I \rho_I m_{I}, \\
\dfrac{d}{dt} \rho_R m_{R} &= \gamma_I \rho_I m_{I}.
\end{split}
\end{equation}
In view of the fact that the equilibrium distribution of the operator \eqref{eq:FPQJ} is either a Gamma distribution, if $\delta = 1$, or an inverse Gamma distribution, if $\delta = -1$, in the limit $\tau \to 0^+$, we have the following equilibrium closure
\begin{equation}\label{eq:closure_L1}
\int_{\R_+} x^{2}f_S(x,t)dx\approx  \int_{\R_+} x^{2}\rho_S(t) f_S^q(x,t)dx = 
\begin{cases}\vspace{0.25cm}
 \dfrac{\lambda+1}{\lambda}\rho_Sm_S^2 & \delta = 1 \\
 \dfrac{\lambda}{\lambda-1} \rho_S m_S^2 & \delta = -1. 
\end{cases}
\end{equation}
being $\lambda = \alpha/\sigma^2$. Hence, we may couple the evolution of mass fractions with the evolution of the first order moment which, in its closed form reads
\begin{equation}\label{eq:mean_L1}
\begin{split}
\dfrac{d}{dt} m_{S} &= -\beta_1 m^2_{S} \rho_I m_{I} (\Lambda - 1),  \\
\dfrac{d}{dt} m_{I} &= \beta_1 \rho_S m_{S} m_{I}(\Lambda m_{S}-m_{I}), \\
\dfrac{d}{dt} m_{R} &= \gamma_I \frac{\rho_I}{\rho_R} (m_{I}-m_{R}), 
\end{split}
\end{equation}
being $\Lambda=((\lambda+\delta)/\lambda)^\delta>1$ for any $\delta = \pm 1$.
The evolution of mass fractions \eqref{eq:mass_L1}, coupled with the closed evolution of first order moments \eqref{eq:mean_L1}, forms a new system of macroscopic equations that inherits the properties of the microscopic agent-based system describing the dynamics of contact formation. The evolution of the mean number of contacts of the susceptible population is always decreasing in time, therefore the maximum number of contacts of the infected population is  
\begin{equation}
\label{eq:maxI_L1}
\bar m_{I} = 
\begin{cases}\vspace{0.25cm}
\dfrac{\lambda+1}{\lambda} m_{S}(0) & \delta = 1, \\
\dfrac{\lambda}{\lambda-1} m_{S}(0) & \delta = -1.
\end{cases}
\end{equation}
Furthermore, we can observe that 
\[
\bar m_{I}^{(\delta = -1)} > \bar m_{I}^{(\delta = 1)},
\]
and the number of contacts of the infected population corresponding to an inverse Gamma closure ($\delta = -1$) are higher compared to the number of contacts obtained with a Gamma closure ($\delta = 1$). 

\begin{figure}\centering
\includegraphics[scale = 0.3]{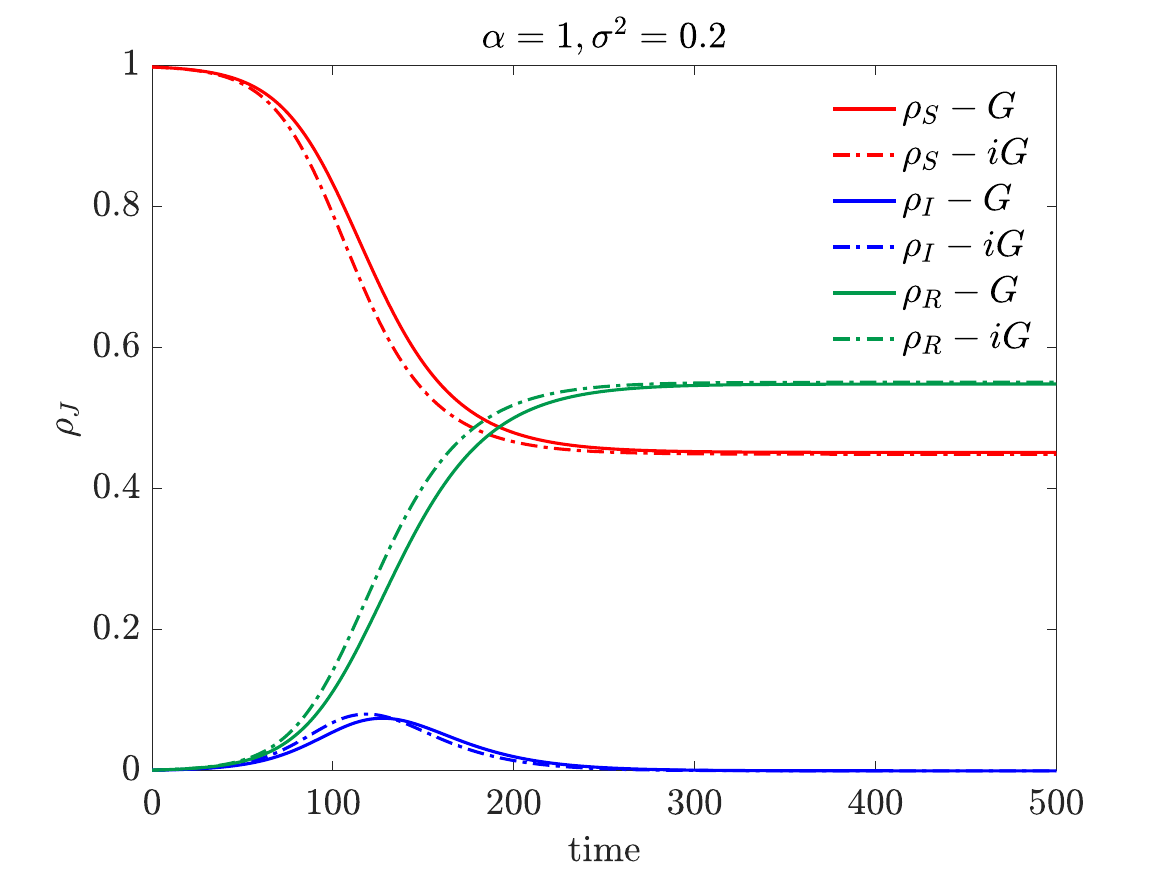}
\includegraphics[scale = 0.3]{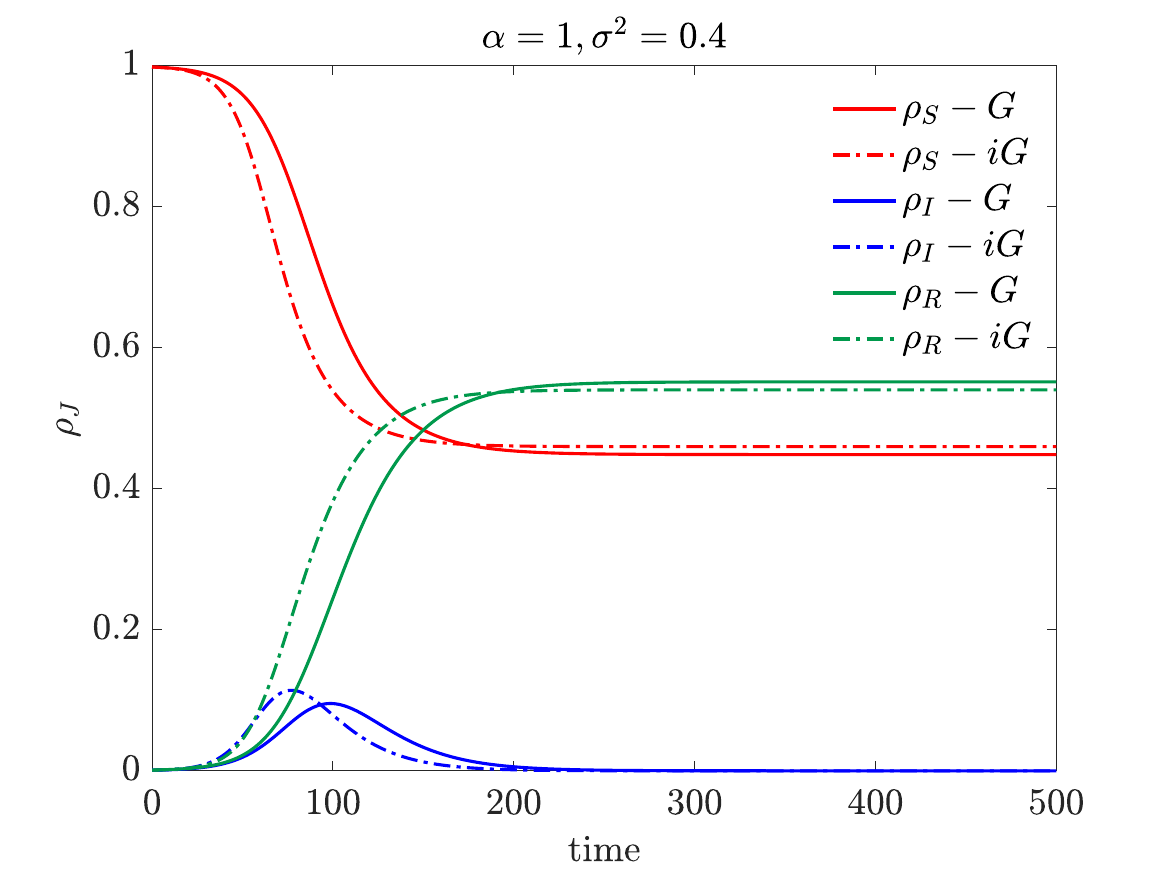}\\
\includegraphics[scale = 0.3]{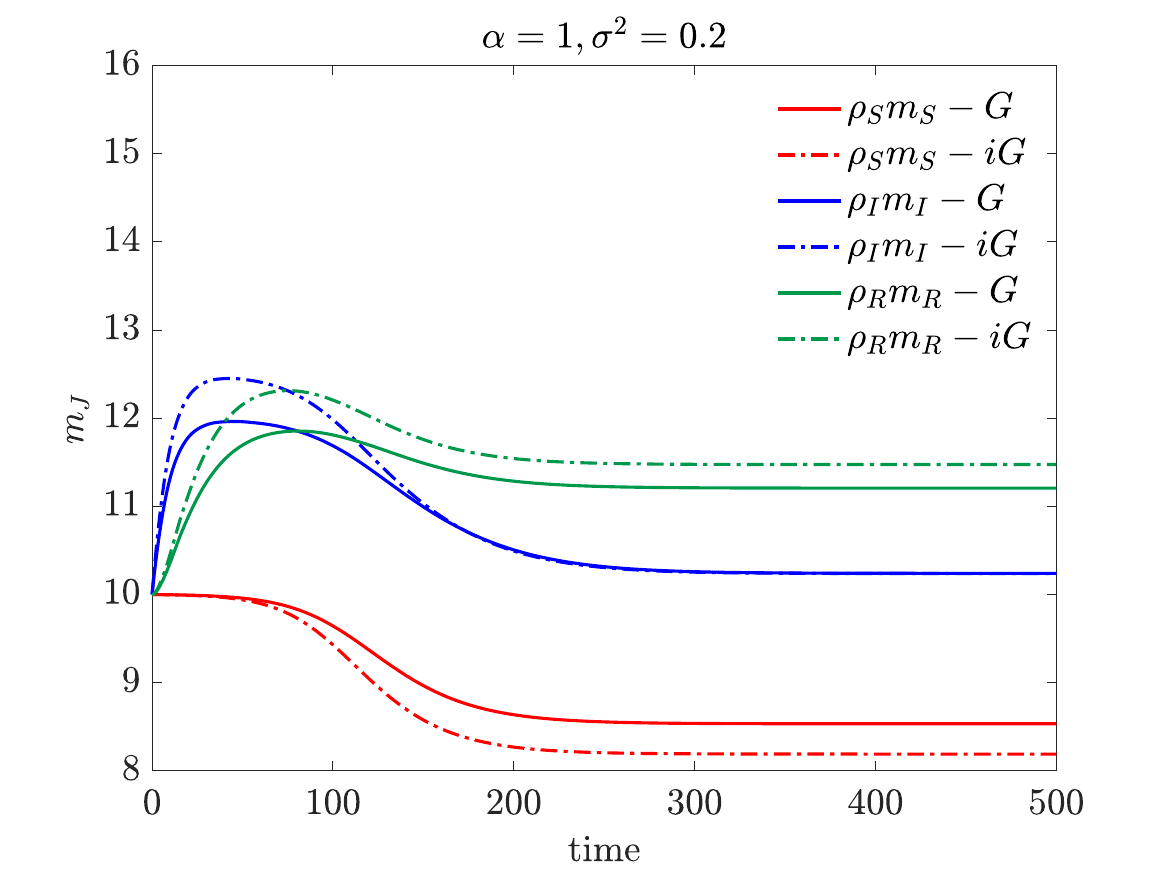}
\includegraphics[scale = 0.3]{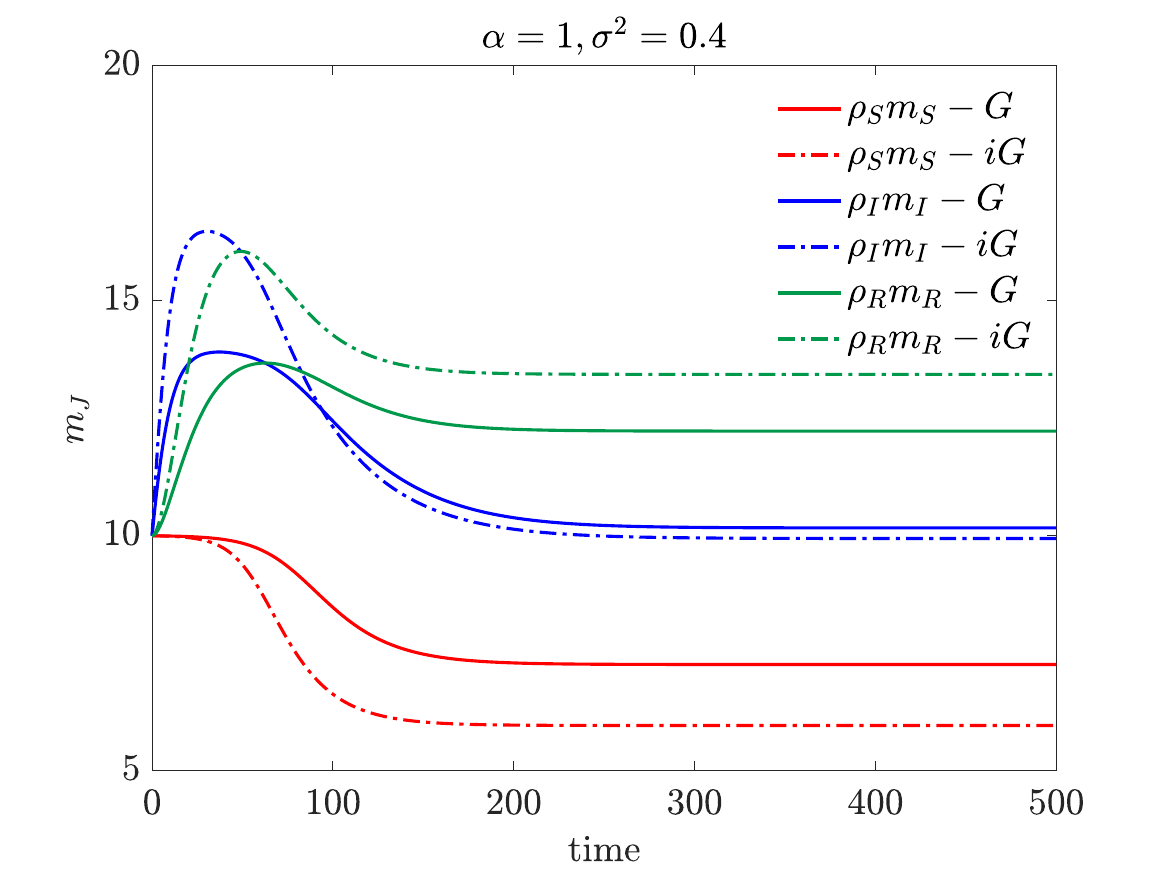}
\caption{\textbf{Case $L=1$}: Epidemic trajectories of the system \eqref{eq:mass_L1}-\eqref{eq:mom_L1} with closure \eqref{eq:closure_L1}. The continuous line corresponds to the evolution of the system obtained with a \textcolor{black}{Gamma closure ($G$)} and the dashdotted line the evolution of the system with \textcolor{black}{inverse Gamma closure ($iG$)}. As initial condition we considered as initial condition $m_J(0) = 10$, $J \in \mathcal C$, and $\rho_I(0) =\rho_R(0) = 10^{-5}$, $\rho_S(0) = 1-\rho_I(0)-\rho_R(0)$ and $\beta_1 = 10^{-3}$, $\gamma_I = 1/14$. From left to right we present the cases $\lambda = 5$ (left), $\lambda = 2.5$ (right).  }
\label{fig:1}
\end{figure}

Consistently with the closure that we obtained, as $\lambda \to +\infty$, which means $\sigma^2 \to 0^+$, the two models describe the same trajectory. It is worth to observe that, if $\sigma\to 0^+$, we get a population with no contact heterogeneities and, therefore, we collapse to the standard SIR model. In other words, the classical SIR model can be obtained through a Dirac delta closure centered on the mean $m_J$. In Figure \ref{fig:1} we depict the epidemic trajectories for several $\sigma^2 =0.2,0.4$ and fixed $\alpha = 1$, in all cases we considered as initial condition $m_J(0) = 10$ and $\rho_I(0) =\rho_R(0) = 10^{-3}$, $\rho_S(0) = 1-\rho_I(0)-\rho_R(0)$. We may observe that for larger $\lambda$ the model obtained with inverse Gamma closure provides a higher mean number of contacts for the infected population. 

Finally, in the case $L = 2$ the evolution of the masses reads
\begin{equation}
\label{eq:mass_L2_}
\begin{split}
\dfrac{d}{dt} \rho_S &= -\beta_1 \rho_S m_{S} \rho_I m_{I} -\beta_2 \rho_S m_{2,S} \rho_I  m_{2,I} \\
\dfrac{d}{dt} \rho_I &=\beta_1 \rho_S m_{S}  m_{I}\rho_I +\beta_2 m_{2,S}\rho_S \rho_I m_{2,I} - \gamma_I \rho_I \\
\dfrac{d}{dt} \rho_R &=\gamma_I \rho_I 
\end{split}
\end{equation}
which can be coupled with the evolution of the mean number of contacts
\begin{equation}
\label{eq:mom_L2}
\begin{split}
\dfrac{d}{dt} \rho_S m_{S} &= -\beta_1 \rho_S m_{2,S} \rho_I m_{I}-\beta_2 \rho_S m_{3,S}\rho_I m_{2,I}  \\
\dfrac{d}{dt} \rho_I m_{I} &= \beta_1 \rho_S m_{2,S} \rho_I m_{I}+\beta_2 \rho_S m_{3,S}\rho_I m_{2,I} - \gamma_I \rho_I  m_{I} \\
\dfrac{d}{dt} \rho_R m_{R} &=\gamma_I \rho_I  m_{I}. 
\end{split}
\end{equation}
As before, in the limit $\tau \to 0^+$ we can obtained a closed evaluation of higher order moments through an equilibrium closure approach, i.e. by exploiting the information in \eqref{eq:gamma} and \eqref{eq:invgamma} for $J \in \{S,I\}$
\begin{equation}
\label{eq:closure_L2}
\begin{split}
\int_{\R} x^2 f_J(x,t)dx \approx  \int_{\R} x^2 f_J^\infty(x,t)dx &= \begin{cases}\vspace{0.25cm}
\dfrac{\lambda+1}{\lambda}\rho_J m_J^2 & \delta = 1 \\
\dfrac{\lambda}{\lambda-1}\rho_J m_J^2 & \delta = -1,
\end{cases} \\
\int_{\R} x^3 f_J(x,t)dx \approx  \int_{\R} x^3 f_J^\infty(x,t)dx &= \begin{cases}\vspace{0.25cm}
\dfrac{(\lambda+1)(\lambda+2)}{\lambda^2}\rho_J m_J^3 & \delta = 1 \\
\dfrac{\lambda^2 }{(\lambda-1)(\lambda-2)}\rho_J m_J^3 & \delta = -1. 
\end{cases}
\end{split}
\end{equation}
For fat-tailed distributions, the obtained closure is derived under the condition that $\lambda > 2$, i.e., $\frac{\alpha}{\sigma^2} > 2$, as the third order moment must be well-defined. Hence, in the limit $\tau \to 0^+$ the evolution of mass fractions reads 
\begin{equation}
\label{eq:mass_L2}
\begin{split}
\dfrac{d}{dt} \rho_S &= -\beta_1 \rho_S m_{S} \rho_I m_{I} -\Lambda^2\beta_2 \rho_S m_{S}^2 \rho_I  m_{I}^2 \\
\dfrac{d}{dt} \rho_I &=\beta_1 \rho_S m_{S}  m_{I}\rho_I +\beta_2\Lambda^2 \rho_S m_{S}^2 \rho_I m_{I}^2 - \gamma_I \rho_I \\
\dfrac{d}{dt} \rho_R &=\gamma_I \rho_I 
\end{split}
\end{equation}
which is coupled to the following system for the evolution of the mean number of contacts
\begin{equation}\label{eq:L2_mean_closed}
\begin{split}
\dfrac{d}{dt} m_{S} &=- m_{S}^2 \rho_I m_{I} \left(  \dfrac{2\beta_1}{2\lambda - 1+\delta} + 2\beta_2 \left(\dfrac{\lambda+\delta}{\lambda}\right)^{2\delta}\dfrac{1}{\lambda-1+\delta} m_{S} m_{I} \right) \\
\dfrac{d}{dt} m_{I} &= \rho_S m_{S}  m_{I} \Bigg[\beta_1 \left(\Lambda m_{S} - m_{I} \right) +  \\
&\quad \beta_2 \Lambda^2 \left(\left( \dfrac{\lambda + 2\delta}{\lambda}\right)^\delta m_{S} - m_{I} \right) m_{S} m_{I} \Bigg], \\
\dfrac{d}{dt}m_{R} &= \gamma_I \dfrac{\rho_I}{\rho_R} \left( m_{I} - m_{R} \right), 
\end{split}
\end{equation}
being $\Lambda = ((\lambda+\delta)/\lambda)^{\delta}>1$ 
with $\delta = \pm 1$. As in the case $L =1$, the mean number of contacts of the susceptible population decreases in time whereas, focussing only on the effects of the second order moments, if $\beta_1=0$ the maximum number of contacts of the infected population is given by 
\begin{equation}
\label{eq:maxI_L2}
\hat m_{I} = 
\begin{cases}\vspace{0.25cm}
\dfrac{\lambda+2}{\lambda}m_{S}(0) & \delta = 1 \\
\dfrac{\lambda}{\lambda-2}m_{S}(0) & \delta =- 1,
\end{cases}
\end{equation}
and we have the following order relation between the obtained maxima in \eqref{eq:maxI_L2}
\[
\hat m_{I}^{(\delta = -1)}> \hat m_{I}^{(\delta = 1)},
\]
for any $\lambda>0$. Hence, the infected population reaches a higher number of contacts in the presence of power-law tails. Furthermore, the obtained maxima for $L = 1$ in \eqref{eq:maxI_L1} and $L = 2$ in \eqref{eq:maxI_L2} are such that 
\[
\bar m_{I}^{(\delta = -1)}< \hat m_{I}^{(\delta = -1)}.
\]
Therefore, the model obtained with taking into account higher order moments produces higher maximum in the mean number of contacts of infected population. 
\begin{figure}\centering
\includegraphics[scale = 0.3]{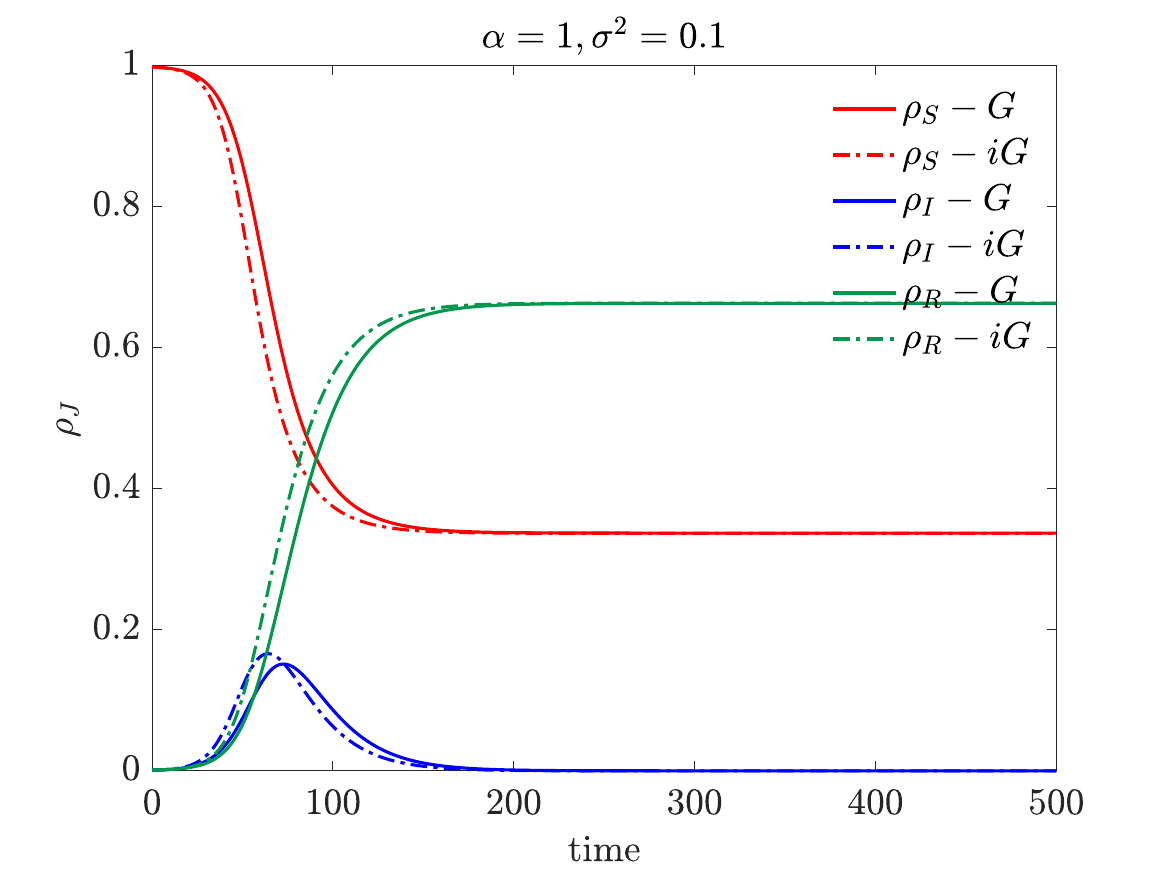}
\includegraphics[scale = 0.3]{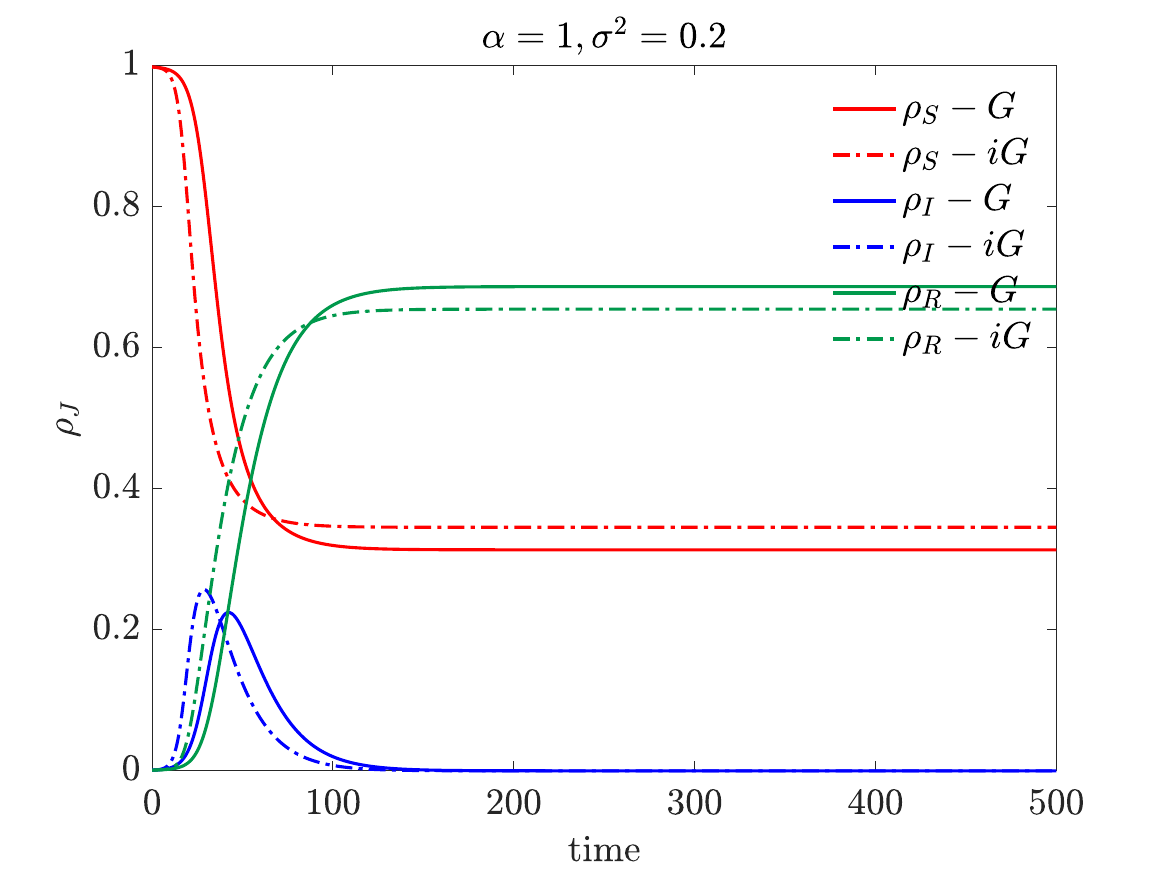}\\
\includegraphics[scale = 0.3]{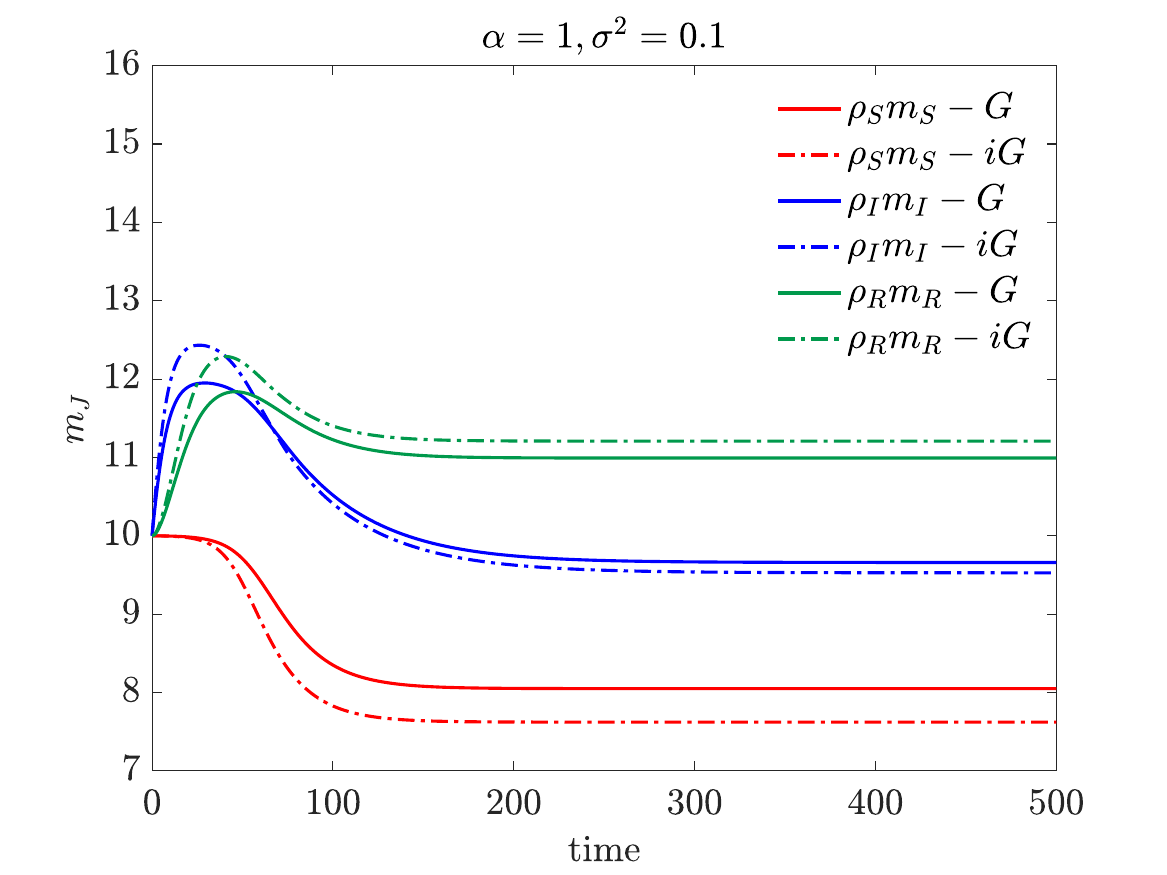}
\includegraphics[scale = 0.3]{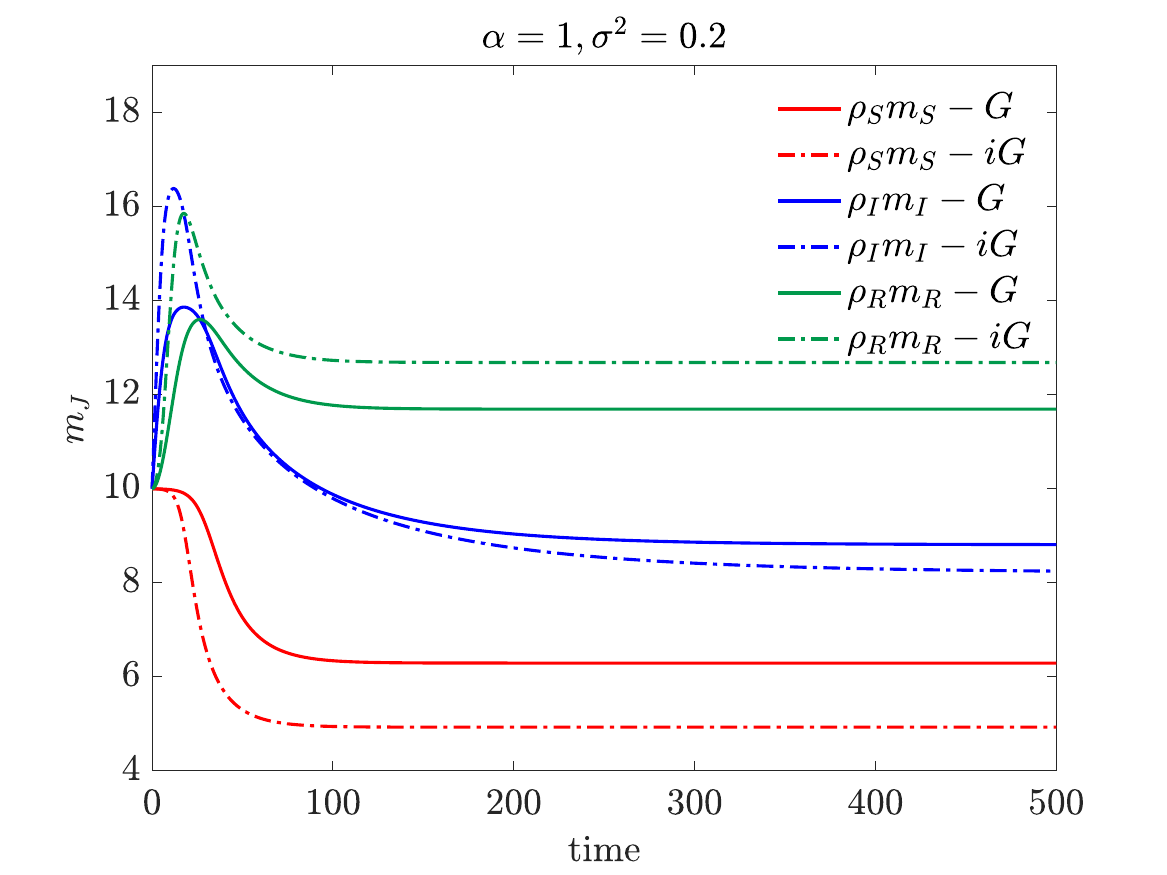}
\caption{\textbf{Case $L=2$}: Epidemic trajectories of the system \eqref{eq:mass_L2}-\eqref{eq:L2_mean_closed} with closure discussed in \eqref{eq:closure_L2}. The continuous line corresponds to the evolution of the system obtained with a \textcolor{black}{Gamma closure ($G$)} and the dashdotted line the evolution of the system with \textcolor{black}{inverse Gamma closure ($iG$)}. As initial condition we considered as initial condition $m_J(0) = 10$, $J \in \mathcal C$, and $\rho_I(0) =\rho_R(0) = 10^{-3}$, $\rho_S(0) = 1-\rho_I(0)-\rho_R(0)$ and $\beta_1 = 0$, $\beta_2 = 10^{-5}$, $\gamma_I = 1/14$. From left to right we present the cases $\lambda = 10$ (left), $\lambda = 5$ (right). }
\label{fig:2}
\end{figure}

In Figure \ref{fig:2} we depict the epidemic trajectories for several $\sigma^2 =0.1,0.2$ and fixed $\alpha = 1$, in all cases we considered as initial condition $m_J(0) = 10$ and $\rho_I(0) =\rho_R(0) = 10^{-3}$, $\rho_S(0) = 1-\rho_I(0)-\rho_R(0)$. We may observe again that for larger $\lambda$ the model obtained with inverse Gamma closure provides higher mean number of contacts for the infected population. 

\section{Control of the tail behaviour}\label{sect:3}

In this section, we will explore the potential to control the dynamics of contact formation in a population, mimicking the role of non-pharmaceutical interventions (NPIs). These interventions aim to reduce the risk factors associated with the transmission of infectious diseases by influencing how individuals interact within the community.

The new kinetic description provides a powerful framework to analyse and quantify the impact of such measures. Specifically, it allows us to highlight the effects of policy-maker actions that directly target the contact distribution within the society. This approach is particularly valuable when only partial information about the population’s contact patterns is available. By leveraging this mathematical approach, we seek to understand how modifying contact dynamics to mitigate the spread of an infection and guide the design of more effective intervention strategies.

In \cite{DTZ,FMZ}, it has been shown that the design and implementation of selective policies can significantly influence the behaviour of contact distributions within a population. These selective policies are measures specifically tailored to act differentially on individuals based on their number of contacts, targeting specific subgroups rather than applying uniform interventions. The introduction of a penalization mechanism inversely proportional to the number of contacts can effectively alter the shape of the contact distribution. This approach enables a transition from a fat-tailed distribution, which is characterized by a high probability to find individuals with large number of contacts, to a slim-tailed distribution, such probability is negligible. This shift in the tail behaviour of the contact distribution has significant implications for public health strategies, as it suggests that selective interventions can play a crucial role in reducing the potential for super-spreading events and enhancing the overall effectiveness of non-pharmaceutical interventions in controlling disease transmission.

To this end, in the following we concentrate on the case $\delta = -1$ which produces inverse Gamma distributions \eqref{eq:invgamma} that are characterised by power-law tails.  To steer the system towards a prescribed number of contacts we proceed as follows: we consider a simple control protocol resulting from the following minimization 
\begin{equation}
\label{eq:cost}
\textrm{argmin}_{u \in \mathcal U} \left[ \dfrac{1}{2}\left\langle(x^\prime-x_T)^2 \right\rangle + \dfrac{\nu}{2}u^2 \right],
\end{equation}
\textcolor{black}{where $\mathcal U$ defines the class of admissible controls such that $x^\prime \ge 0$}, \textcolor{black}{$x_T\in\R^+$ is the target contact number, and $\nu>0$ is a penalisation coefficient}. \textcolor{black}{In \eqref{eq:cost}, the cost increases quadratically with the distance to a target state. This choice reflects the fact that greater effort is required to contain populations with a large average number of contacts. The cost could be different than zero, allowing the existence of agents with a controlled non-zero number of contacts.} The minimisation \eqref{eq:cost} is subject to the two possible dynamics that define the post-interaction state $x^\prime\ge0$. Hence, the minimisation of the cost in \eqref{eq:cost} can be obtained through a Lagrange multiplier approach, see \cite{APTZ,PTZ}. In more detail, in the following we will consider the following transition laws:
\begin{itemize}
\item[$A)$] Additively controlled transition 
\begin{equation}\label{eq:additive}
x^\prime = x  + \epsilon( \Psi(x/m_J)x + u) + x \eta,
\end{equation}
\textcolor{black}{with $\Psi(x/m_J)$ given by \eqref{eq:psi}.} Such class of controls have been classically considered in several constrained multiagent dynamics to enforce the emergence of prescribed patterns through the addition of external force in the dynamics, see e.g. \cite{ACFK,APTZ} and the references therein. In this case, the minimisation of \eqref{eq:cost} gives
\[
u^* = - \frac{\epsilon}{\nu + \epsilon^2} \left((x-x_T) + \epsilon\Psi(x/m_J) x \right)
\]
which leads to the constrained interaction  
\begin{equation}
\label{eq:controlledA}
x^\prime = x + \dfrac{\nu \epsilon }{\nu + \epsilon^2 }\Psi(x/m_J)x + \dfrac{\epsilon^2 }{\nu+ \epsilon^2}(x_T-x)+x \eta.
\end{equation}

\item[$B)$] Interaction-driven controlled transition
\begin{equation}\label{eq:multiplicative}
x^\prime = x + \epsilon \Psi(x/m_J)x u + x\eta,
\end{equation}
which tries to modify the transition protocol of each single agent to steer the mean number of contacts towards $x_T>0$. \textcolor{black}{Again, the function $\Psi(x/m_J)$ in \eqref{eq:multiplicative} is the one defined in \eqref{eq:psi}.} Hence, the minimisation of \eqref{eq:cost} gives 
\[
u^* =  -\frac{ \epsilon \Psi(x/m_J)x }{\nu + \epsilon^2(x\Psi(x/m_J))^2} (x-x_T).
\]
and the constrained interaction is given by 
\begin{equation}\label{eq:controlledB}
x^\prime = x - \frac{ \epsilon^2(x\Psi(x/m_J))^2 }{\nu + \epsilon^2(x\Psi(x/m_J))^2} (x-x_T)  +x \eta.
\end{equation}
This control has been introduced in \cite{PTZ} in a related modelling setting.
\end{itemize}
It is worth to notice that in both cases \eqref{eq:additive}-\eqref{eq:multiplicative}, in the limit $\nu\to 0^+$ we obtain an instantaneous control of the transition which is driven towards $x_T>0$. \textcolor{black}{Anyway, the introduced controlled interactions have a different impact on the moments of the emerging equilibrium distribution at the kinetic level.}
\subsection{Fokker-Planck model and controlled equilibria}\label{sect:31}
\textcolor{black}{We introduce the notation $f_J^{(A)}(x,t)$, $f_J^{(B)}(x,t)$ to denote the distribution of the number of contacts of a multiagent system characterised by microscopic constrained interaction expressed by the controls $(A)$ and $(B)$, respectively defined in \eqref{eq:additive} and \eqref{eq:multiplicative}, in the compartment $J$. We fix $\delta=-1$ from now on when considering the controlled scenario, since we are interested in the control of fat-tailed distributions. Proceeding as in Section \ref{sect:contact}, we obtain the evolution of observable quantities through the following space homogeneous Boltzmann-type equation
}
\begin{equation}
\label{eq:boltzmann_control}
\textcolor{black}{\dfrac{d}{dt} \int_{\mathbb R_+}\varphi(x)f_J^{(H)}(x,t)dx = \left\langle \int_{\mathbb R_+}(\varphi(x^\prime)-\varphi(x))f_J^{(H)}(x,t)  dx\right\rangle,}
\end{equation}
where $ H = A,B$, \textcolor{black}{and we note that the kernel \eqref{eq:kernel} is unitary for $\delta=-1$}.  By introducing the scaling of the penalization term $\nu \to \epsilon \nu$, in the quasi-invariant limit and for $\epsilon \to 0^+$,   as observed in \cite{PTZ}, we can obtain a Fokker-Planck-type equation with a modified drift term that takes into account the presence of the control. In terms of the controlled density \textcolor{black}{$f_J^{(A)}(x,t)$}, this equation reads 
\begin{equation} \label{eq:FPQJA}
\begin{split}
&\textcolor{black}{\partial_t f_J^{(A)}(x,t) =} \\
&\textcolor{black}{ \quad\partial_x  \left\{\left[ \dfrac{\alpha}{2\delta}\left( \dfrac{m_J}{x} -1\right)x + \dfrac{x-x_T}{\nu}  \right ]f_J^{(A)}(x,t) + \dfrac{\sigma^2}{2}\partial_x(x^{2}f_J^{(A)}(x,t))\right\} = }\\
&\quad \textcolor{black}{Q^{(A)}(f_J^{(A)})(x,t)}
\end{split}
\end{equation}
which has the steady state
\begin{equation}
\label{eq:steadyA}
\textcolor{black}{f_J^{(A),\infty}(x)= C_{-1,\sigma^2,\alpha} x^{\frac{\lambda}{ \delta}-2-\frac{2}{\sigma^2 \nu}} \exp\left\{ - \lambda \dfrac{m_J}{x} - \dfrac{2}{\sigma^2 \nu}\dfrac{x_T}{x}\right\},}
\end{equation}
where the constant \textcolor{black}{$C_{-1,\sigma^2,\alpha}>0$} is chosen such as the mass of the density is unitary. 
\textcolor{black}{We observe that the equilibrium distribution $f_J^{(A),\infty}$ is still a fat-tailed distribution and, therefore, the action of the control is not capable of modifying the tails of the distribution.}

On the other hand, in terms of the density \textcolor{black}{$f_J^{(B)}(x,t)$} the Fokker-Planck equation derived from \eqref{eq:boltzmann_control} reads
\begin{equation} \label{eq:FPQJB}
\begin{split}
&\textcolor{black}{\partial_t f_J^{(B)}(x,t) = }\\
&\textcolor{black}{\quad\partial_x  \left\{\dfrac{x^{2}}{\nu}\left[ - \dfrac{\alpha}{2}\left(\dfrac{m_J}{x} -1\right)\right ]^2(x-x_T) f_J^{(B)}(x,t) + \dfrac{\sigma^2}{2}\partial_x(x^{2}f_J^{(B)}(x,t))\right\}=} \\
&\quad \textcolor{black}{Q^{(B)}(f_J^{(B)})(x,t)}
\end{split}
\end{equation}
\textcolor{black}{whose equilibrium distribution, in the case $\delta \to -1$, reads}
\[
\textcolor{black}{f_J^{(B),\infty}(x) = C_{-1,\sigma^2,\alpha}x^{-2-\ell} \exp\left\{-\dfrac{\alpha^2}{2\sigma^2\nu} \left[ - (2m_J+x_T)x + \dfrac{x^2}{2}+\dfrac{m_J^2x_T}{x} \right]\right\},}
\]
where $\ell = \frac{\alpha^2}{2\sigma^2\nu}(m_J^2 + 2x_T m_J)$. The obtained equilibrium density can be equivalently rewritten as follows
\begin{equation}
\label{eq:steadyBB}
\begin{split}
f_J^{(B),\infty}(x) =  C_{-1,\sigma^2,\alpha}x^{-2-\ell} \mathcal N\left(2m_J + x_T,\dfrac{2\sigma^2\kappa}{\alpha^2}\right) \chi(x \ge0)\times \\
\exp\left\{-\dfrac{\alpha^2}{2\sigma^2\nu} \left[-\dfrac{(2m_J^2+x_T)^2}{2}+\dfrac{m_J^2 x_T}{x} \right] \right\},
\end{split}
\end{equation}
being $\chi(\cdot)$ the indicator function and $\mathcal N(\cdot,\cdot)$ a Gaussian distribution with mean defined in the first entry and variance in the second entry. Hence, in the case of the control $(B)$ the equilibrium distribution exhibits slim tails. 

\begin{remark}\label{rem:control}
We may notice that the evolution of the controlled densities \textcolor{black}{$f_J^{(H)}(x,t)$, $ H = A,B$}, is obtained in terms of the Fokker-Planck equations \eqref{eq:FPQJA}-\eqref{eq:FPQJB} which conserve only the mass. Hence, the momentum is not a conserved quantity for the presented control strategies. 
\end{remark}

\subsection{Controlled macroscopic equations}\label{sect:32}
In this section we seek to obtain a closed system for the evolution of the main moments of the controlled kinetic system
\begin{equation}
\label{eq:kinetic_control}
    \begin{cases}
\partial_t f^{(H)}_S(x,t) = -K(f^{(H)}_S,f^{(H)}_I)(x,t) + \dfrac{1}{\tau}Q^{(H)}_S(f^{(H)}_S)(x,t), \\
\partial_t f^{(H)}_I(x,t) =K(f^{(H)}_S,f^{(H)}_I)(x,t)-\gamma_I f^{(H)}_I(x,t) + \dfrac{1}{\tau}Q^{(H)}_I(f^{(H)}_I)(x,t) \\
\partial_t f^{(H)}_R(x,t) = \gamma_I f^{(H)}_I(x,t) +  \dfrac{1}{\tau}Q^{(H)}_R(f^{(H)}_R)(x,t),
    \end{cases}
\end{equation}
with $H = A,B$ and $Q^{(H)}(f^{(H)})(x,t)$ defined in \eqref{eq:FPQJA} and \eqref{eq:FPQJB}. As observed in Remark \ref{rem:control}, the only conserved quantity of the operators $Q^{(H)}(f^{(H)})$ is the mass. 

To determine the evolution of mass fractions we integrate \eqref{eq:kinetic_control} to obtain the evolution of the controlled macroscopic equations. In the case $L = 1$ we get
\begin{equation}
\label{eq:macrocontrol_L1}
\begin{split}
\dfrac{d}{dt} \rho^{(H)}_S &= -\beta_1 \rho^{(H)}_S m^{(H)}_{S} \rho^{(H)}_I m^{(H)}_{I}  \\
\dfrac{d}{dt} \rho^{(H)}_I &=\beta_1 \rho^{(H)}_S m^{(H)}_{S}  m^{(H)}_{I}\rho^{(H)}_I  - \gamma_I \rho^{(H)}_I \\
\dfrac{d}{dt} \rho^{(H)}_R &=\gamma_I \rho^{(H)}_I 
\end{split}
\end{equation}
where, in the limit $\tau \to 0^+$, we recover the information on the constrained mean value as follows
\[
m^{(H)}_{1,J}\rho^{(H)}_J \approx \rho_J^{(H)}(t)\int_{\R^+}xf^{(H),\infty}_{J}(x)dx. 
\]
being $f^{(H),\infty}(x,t)$ the steady state of the controlled problem which depends on the adopted control strategy, see \eqref{eq:steadyA}-\eqref{eq:steadyBB}. Therefore, the case $L = 1$ is not sufficient to observe the behaviour of the tail distribution of the controlled contact formation dynamics. 

On the other hand, in the case $L = 2$ we get
\begin{equation}
\label{eq:macrocontrol_L2}
\begin{split}
\dfrac{d}{dt} \rho^{(H)}_S &= -\beta_1 \rho^{(H)}_S m^{(H)}_{S} \rho^{(H)}_I m^{(H)}_{I} -\beta_2 m^{(H)}_{2,S}\rho^{(H)}_S  m^{(H)}_{2,I}\rho^{(H)}_I \\
\dfrac{d}{dt} \rho^{(H)}_I &=\beta_1 \rho^{(H)}_S m^{(H)}_{S}  m^{(H)}_{I}\rho^{(H)}_I +\beta_2 m^{(H)}_{2,S}\rho^{(H)}_S \rho^{(H)}_I m^{(H)}_{2,I} - \gamma_I \rho^{(H)}_I \\
\dfrac{d}{dt} \rho^{(H)}_R &=\gamma_I \rho^{(H)}_I 
\end{split}
\end{equation}
where, as before, in the limit $\tau\to 0^+$, we can close the system through an equilibrium closure approach
\begin{equation}
    \label{eq:closure_controlL2}
\begin{split}
m^{(H)}_{J}\rho^{(H)}_J &\approx \rho_J^{(H)}(t)\int_{\R^+}xf^{(H),\infty}_{J}(x)dx \\
m^{(H)}_{2,J}\rho^{(H)}_J &\approx \rho_J^{(H)}(t)\int_{\R^+}x^{2}f^{(H),\infty}_{J}(x)dx. 
\end{split}
\end{equation}

We observe that the control strategy $(B)$, which ensures the contact distribution $f^{(B)}(x)$ in \eqref{eq:steadyBB} is slim-tailed, consistently results in a finite second-order moment for the emerging controlled distribution. In contrast, control strategy $(A)$ fails to guarantee the finiteness of the second-order moment for the corresponding controlled distribution, $f^{(A),\infty}(x)$ in \eqref{eq:steadyA}. \textcolor{black}{At the  level of epidemic dynamics, in the presence of large multiagent systems,} this indicates that interaction-driven controls are more effective than uniform-type controls in reducing the probability of finding agents with high number of contacts. 

\section{Numerical tests}\label{sect:4}
In this section, we present several numerical results. 
In particular, in Test 1 we check for the long-time accordance between the numerical resolution of the Boltzmann-type equation \eqref{eq:boltz} with binary interaction rules given by \eqref{eq:micro}-\eqref{eq:controlledA}-\eqref{eq:controlledB} in the quasi invariant limit, and the limiting Fokker-Planck equations \eqref{eq:FPQJ}-\eqref{eq:FPQJA}-\eqref{eq:FPQJB}. Numerical results, in both the uncontrolled and controlled scenario, are compared with the corresponding analytical equilibrium distributions.
In Test 2 we investigate the behaviour of the tails of the contact distribution in the absence of epidemic exchange. We study the dependence of the different control strategies with respect to the penalization $\nu$ and the parameter $\lambda$. 
In Test 3 we check the consistency of the macroscopic limit, namely the accordance between the system \eqref{eq:kinetic} in the limit $\tau\to0$ with $L=2$, and the system \eqref{eq:mass_L2}-\eqref{eq:L2_mean_closed} obtained with the moment closure described in Section \ref{sect:23}.
In Test 4 we numerically solve the system \eqref{eq:kinetic_control} for $H=A,B$ and different penalization and scale parameters. We compare the results with the uncontrolled scenario to study the effectiveness of the control strategies to reduce the spread of the epidemics. 

The Boltzmann-type equation \eqref{eq:boltz} with binary interaction rules given by \eqref{eq:micro}-\eqref{eq:controlledA}-\eqref{eq:controlledB} is solved with a classical DSMC scheme \cite{PT}. \textcolor{black}{To construct the scheme, we discretize the time domain with a step size $\Delta t>0$ such that $t^n=n\Delta t$, and we denote by $f^n_J(x)$ an approximation of $f_J(x,t^n)$ at the $n$-th time step. We then substitute the kernel with $B_\Sigma(x)=\min\{B(x),\Sigma\}$, where $\Sigma>0$ is an upper bound for the kernel, and we rearrange the collisional operator to highlight the gain and loss part. Applying a forward Euler method to approximate the time derivative, we have
\begin{equation}\label{eq:DSMC}
f^{n+1}_J = \left( 1 - \frac{ \Sigma \Delta t }{\epsilon} \right) f^n_J + \frac{ \Sigma \Delta t }{\epsilon} \frac{P(f^n)}{\Sigma},
\end{equation}
where $P(\cdot)$ is 
\[
P(f)(x,t)=\left\langle\int_{\R^+} B_\Sigma(x)  \dfrac{1}{{}^\prime \mathcal{J}} f_J({}^\prime x,t) dx \right\rangle + (1-\Sigma) f_J(x).
\] 
Provided $\Delta t<\epsilon/\Sigma$, \eqref{eq:DSMC} is a convex combination of probability density function, and we can apply the classical DSMC approach. For further details, we refer to \cite{PT}}

For the Fokker-Planck equations \eqref{eq:micro}-\eqref{eq:controlledA}-\eqref{eq:controlledB} in the absence of epidemic exchange, \textcolor{black}{i.e. Test 1-2}, we adopt a structure-preserving (SP) implicit scheme, as proposed in \cite{PZ}. To numerically solve the system \textcolor{black}{for Test 3-4}, we rely on a time splitting technique \cite{FMZ}. In particular, we first rewrite the system in vector form
\[
\frac{\partial \textbf{f}^{(H)}}{\partial t}(x,t) = \textbf{P}(x,\textbf{f}^{(H)}(x,t)) + \frac{1}{\tau} \textbf{Q}^{(H)}(\textbf{f}^{(H)}(x,t)),
\]
where $\textbf{f}^{(H)}=\{f^{(H)}_J\}_J$ and $\textbf{Q}^{(H)}=\{Q^{(H)}_J\}_J$, with $J=\{S,I,R\}$ and $H=A,B$, and $\textbf{P}$ is the vector representing the mass exchange between the compartments. 
Denoting by $\textbf{f}^{(H),n}(x)$ an approximation of $\textbf{f}^{(H)}(x,t^n)$ at the $n$-th time step, the splitting reads
\begin{equation} \label{eq:FP_step}
\textrm{Fokker-Planck contact dynamics: } 
\begin{cases}
    \dfrac{\partial\textbf{f}^{(H),*}}{\partial t} = \frac{1}{\tau} \textbf{Q}^{(H)}(\textbf{f}^{(H),*}) \\ \\
    \textbf{f}^{(H),*}(x,0)=\textbf{f}^{(H),n}(x)
\end{cases}
\end{equation}
\begin{equation} \label{eq:epidemic_dynamics}
\textrm{Epidemic dynamics: }
\begin{cases}
    \dfrac{\partial\textbf{f}^{(H),**}}{\partial t} = \textbf{P}(x,\textbf{f}^{(H),**}) \\ \\
    \textbf{f}^{(H),**}(x,0)=\textbf{f}^{(H),*}(x,\Delta t).
\end{cases}
\end{equation}
Finally, the solution at time $t^{n+1}$ is $f^{(H),n+1}(x)=f^{(H),**}(x,\Delta t)$. The Fokker-Planck contact dynamics \eqref{eq:FP_step} is solved with an implicit SP approach. The system of epidemic exchange \eqref{eq:epidemic_dynamics} is solved with a fourth order Runge-Kutta method.

\subsection{Test 1: Boltzmann to Fokker-Planck}
\begin{figure}[t!]
\centering
\includegraphics[width=0.475\textwidth]{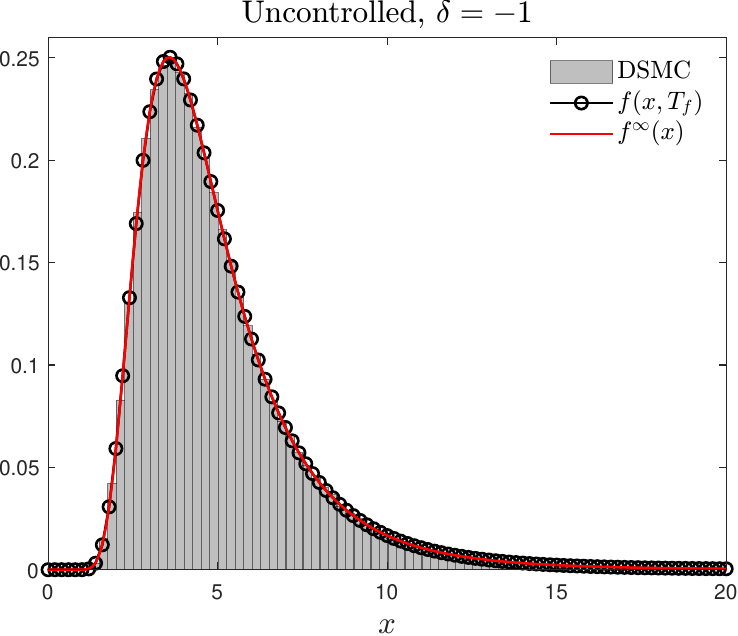}
\includegraphics[width=0.475\textwidth]{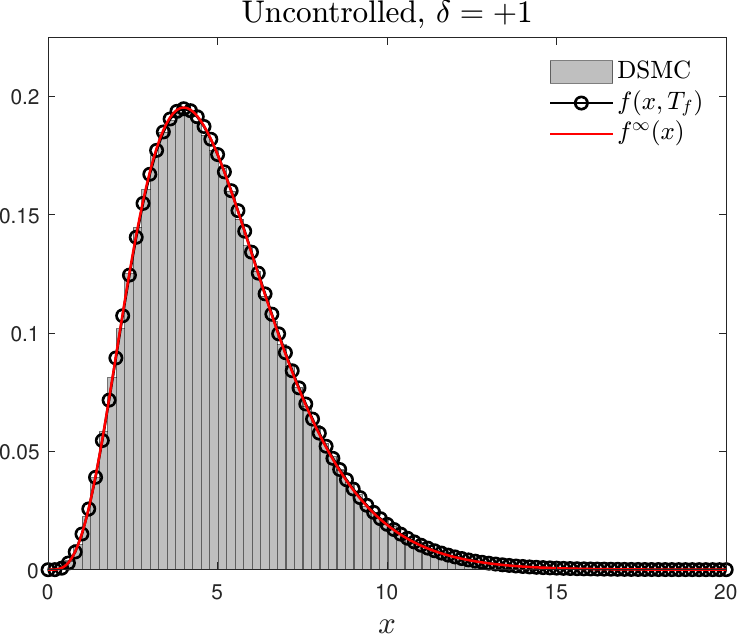} \\
\vspace{5pt}
\includegraphics[width=0.475\textwidth]{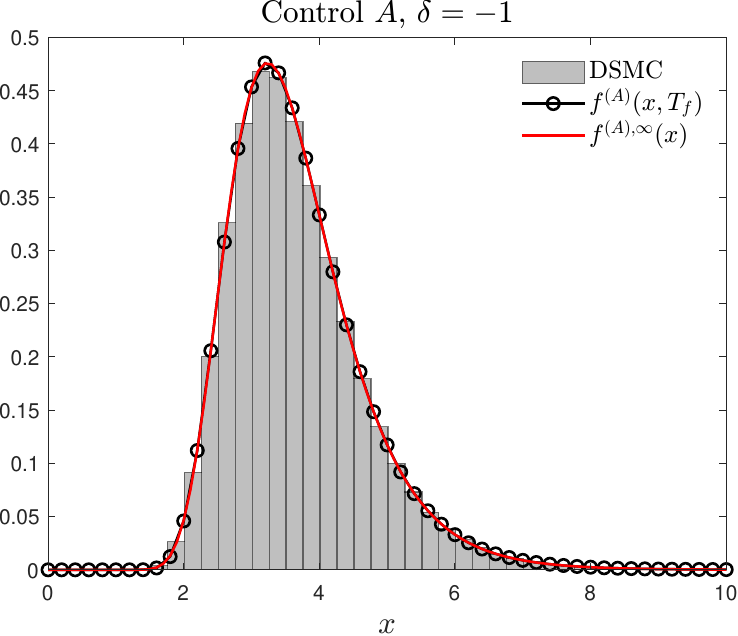}
\includegraphics[width=0.475\textwidth]{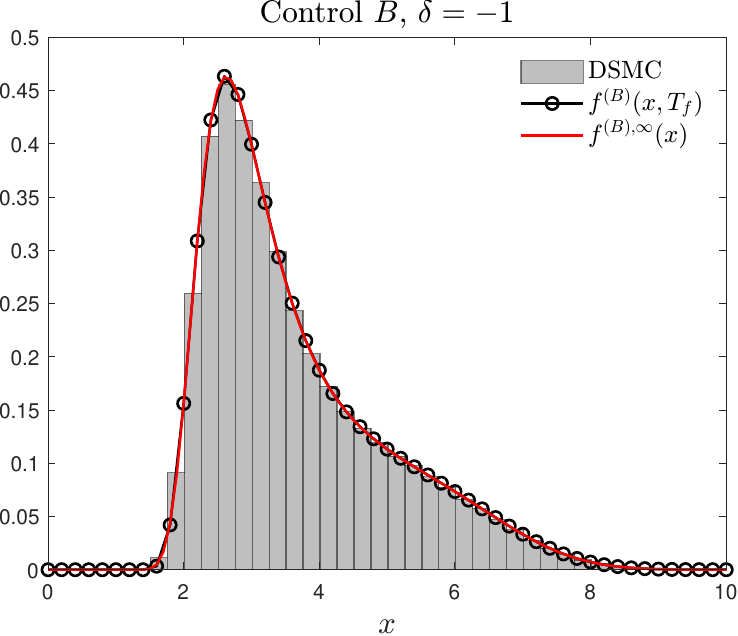}
\caption{Comparison of the long-time behaviour of the numerical solution to the Boltzmann equations (grey bricks) in the quasi in variant interaction limit, the limiting Fokker-Planck equations (black circled lines), and the corresponding equilibrium distributions (solid red lines). Upper row: uncontrolled scenario with $\delta=-1$ (left), and $\delta=+1$ (right), with $T_f=50$ in both cases. Bottom row: controlled scenario with strategies $A$ (left) and $B$ (right), for $\delta=-1$ and with $T_f=20$ in both cases. The $x$-domain is $[0,100]$, and it is cut for visualization purposes. The other parameters are $\alpha=1$, $\sigma^2=0.2$, $m_J=5$, $x_T=3$, $\nu=1$, $\Delta t=\epsilon=0.01$ for $\delta=-1$, and $\Delta t = \epsilon/10=0.001$ for $\delta=+1$.}
\label{fig:test1_BtoFP}
\end{figure}
In this test, \textcolor{black}{we consider the dynamics of a single Boltzmann-type equation, in the quasi invariant collision limit, to check the consistency of the large time solution with the steady state of the derived Fokker-Planck equation,} \textcolor{black}{in both the uncontrolled and controlled scenarios. In more detail, we first compare equation \eqref{eq:boltz} expressing the microscopic transition  \eqref{eq:micro}, with the large time solution of  \eqref{eq:FPQJ},} for both $\delta=+1,-1$. Then, we consider the controlled scenario with $\delta=-1$ and we solve \eqref{eq:controlledA}-\eqref{eq:controlledB} and \eqref{eq:FPQJA}-\eqref{eq:FPQJB}, respectively. All the numerical results are also compared with the corresponding analytical equilibrium distribution.

We fix the parameters $\alpha=1$, $\sigma^2=0.2$, $m_J=5$, the $x$-domain $[0,100]$, and $x_T=3$, $\nu=1$ for the controlled cases. For Maxwellian molecules corresponding to the case with $\delta=-1$, we choose $\Delta t=\epsilon=0.01$, while for the non-Maxwellian case, i.e. $\delta=+1$, we fix $\epsilon=0.01$ and we consider an upper bound for the kernel $B(\cdot)$ such that $\Delta t = \epsilon/10=0.001$. The time domain is $[0,T_f]$, with $T_f=50$ in the uncontrolled case, and $T_f=20$ in the controlled scenario. For the DSMC simulations, we consider $N=10^6$ particles and a histogram reconstruction with $400$ bins in the $x$-domain. For the resolution of the Fokker-Planck equation, we choose $\Delta x=0.1$. The initial condition, in all the simulations, is a normalized uniform distribution $f(x,0)=\mathcal{U}([6,8])$.

In Figure \ref{fig:test1_BtoFP} we may observe the accordance between the numerical solutions to the Boltzmann equations (grey rectangles) in the quasi in variant interaction limit, the limiting Fokker-Planck equations (black circled lines), and the corresponding equilibrium distributions (solid red lines). \textcolor{black}{This test shows that the long-time behaviour of the Boltzmann equation in the quasi invariant limit and of the limiting Fokker-Planck equation are in agreement with the analytical equilibrium distributions, in both the uncontrolled and controlled scenario.}

\subsection{Test 2: Tails behaviour of the contact distribution}
\begin{figure}[t!] 
\centering
\includegraphics[width=0.45\textwidth]{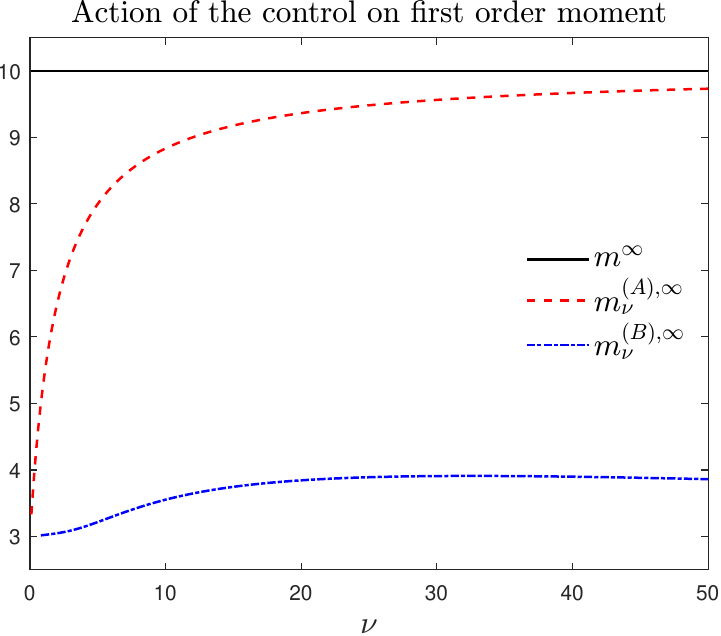}
\includegraphics[width=0.45\textwidth]{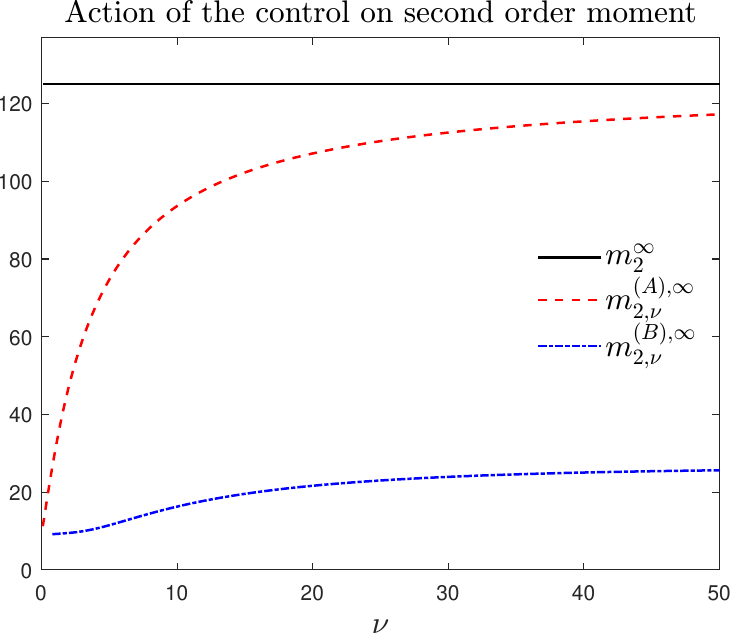}
\caption{First (left) and second (right) order moment at the equilibrium versus the penalization coefficient $\nu$. Solid black lines: uncontrolled scenario independent from the penalization coefficient, thus constant; dashed red lines: additive control strategy $(A)$; dashdotted black lines: multiplicative control strategy $(B)$. The other coefficients are: $\delta=-1$, $\sigma^2=0.2$, $m=10$, $x_T=3$, and $\lambda=3$.}
\label{fig:test2-1}
\end{figure}
In this subsection, we consider a single Fokker-Planck-type equation that describes the formation of the contact dynamics without epidemic exchange. We are interested in understand the tails behaviour in the presence of the control strategies $(A)$ and $(B)$, with respect to the penalization coefficient $\nu$ and the parameter $\lambda=\alpha/\sigma^2$. We concentrate on the case $\delta=-1$, since we have fat-tailed distribution at the equilibrium in the absence of the control. In all the tests, we fix the parameters $\delta=-1$, $\sigma^2=0.2$, $m=10$, and $x_T=3$. 

\begin{figure}[t!] 
\centering
\includegraphics[width=0.475\textwidth]{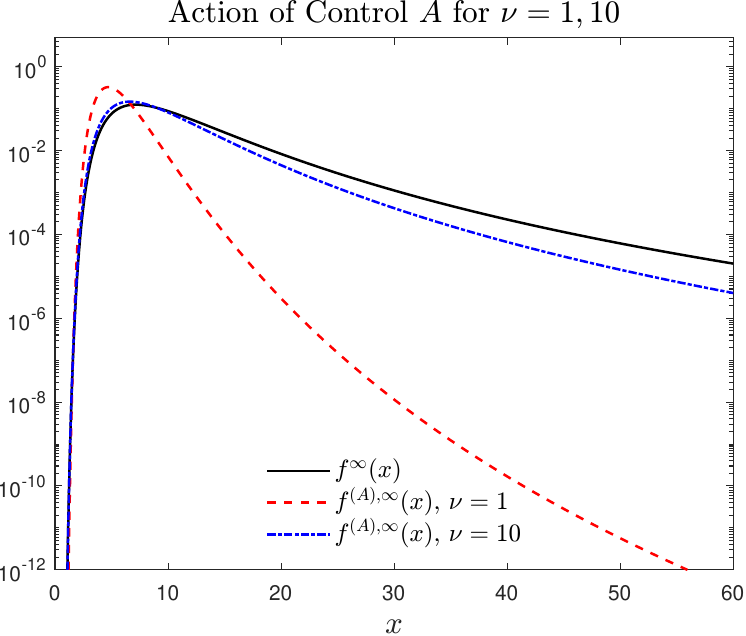}
\includegraphics[width=0.475\textwidth]{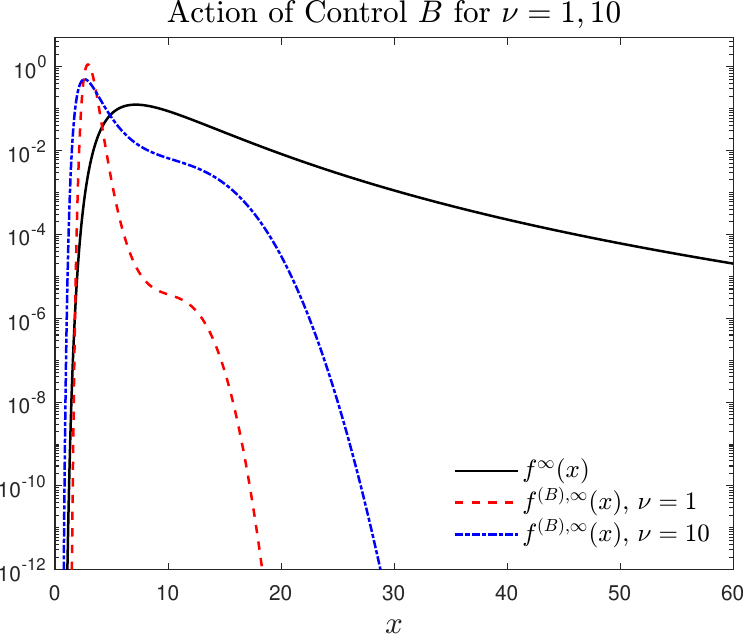}
\caption{Contact distribution at the equilibrium for different penalization coefficients and control strategies. From left to right: $\nu=1,10$. In every plot the solid black line corresponds to the uncontrolled scenario $f^\infty(x)$, the dashed red lines to the additive control strategy $A$, $f^{(A),\infty}(x)$, and the dashdotted black lines to the multiplicative control strategy $B$, $f^{(B),\infty}(x)$. The other coefficients are: $\delta=-1$, $\sigma^2=0.2$, $m=10$, $x_T=3$, and $\lambda=2$.}
\label{fig:test2-2}
\end{figure}

In Figures \ref{fig:test2-1}-\ref{fig:test2-2}, we also fix $\alpha=0.4$ in a way that $\lambda=2$, and we vary the penalization coefficient $\nu$. In particular, in Figure \ref{fig:test2-1} we show the first order moment $m^{\infty}$ and $m^{(H),\infty}_\nu$ and the second order moment $m_2^{\infty}$ and $m^{(H),\infty}_{2,\nu}$ at the equilibrium as a function of $\nu$ for the control strategies $(A)$ and $(B)$. Obviously, in the uncontrolled scenario there is no parametric dependence on $\nu$ and therefore the mean and the energy are constants. We may notice that, as $\nu$ decreases, for both the control strategies the mean approaches the same value, which is the selected target $x_T=3$. We can also observe that the control $(B)$ reduces at the equilibrium both the mean $m^{(H),\infty}_\nu$ and the energy $m^{(H),\infty}_{2,\nu}$ more than the control $(A)$. 
In Figure \ref{fig:test2-2}, we observe the differences in the control of the tails of the distribution at the equilibrium $f^{(H),\infty}(x)$. In fact, as noticed in the previous sections, the control strategy $(B)$ ensures that $f^{(B),\infty}(x)$ is a slim-tailed distribution, while the strategy $(A)$ not. In these plots, we select $\nu=1,10$, and we look at the equilibrium distribution in semilogarithmic scale to highlight the behaviour of the tails.

\begin{figure}[t!] 
\centering
\includegraphics[width=0.475\textwidth]{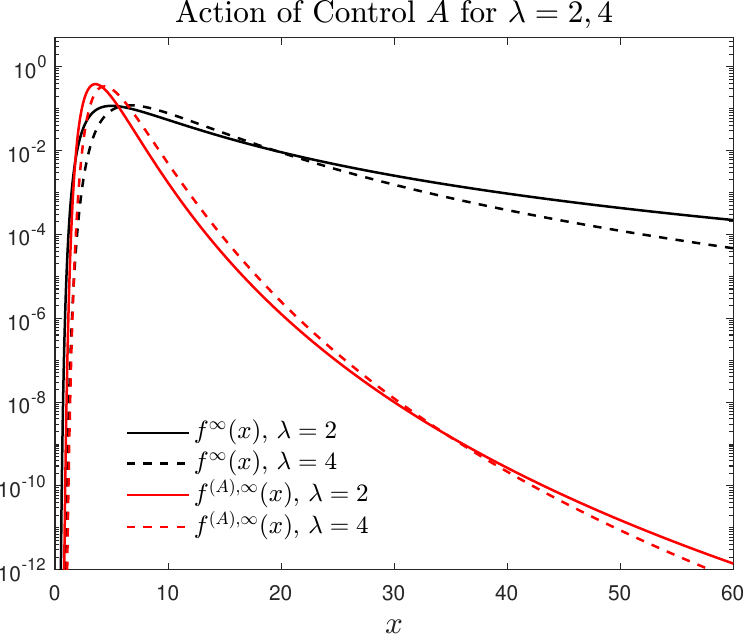}
\includegraphics[width=0.475\textwidth]{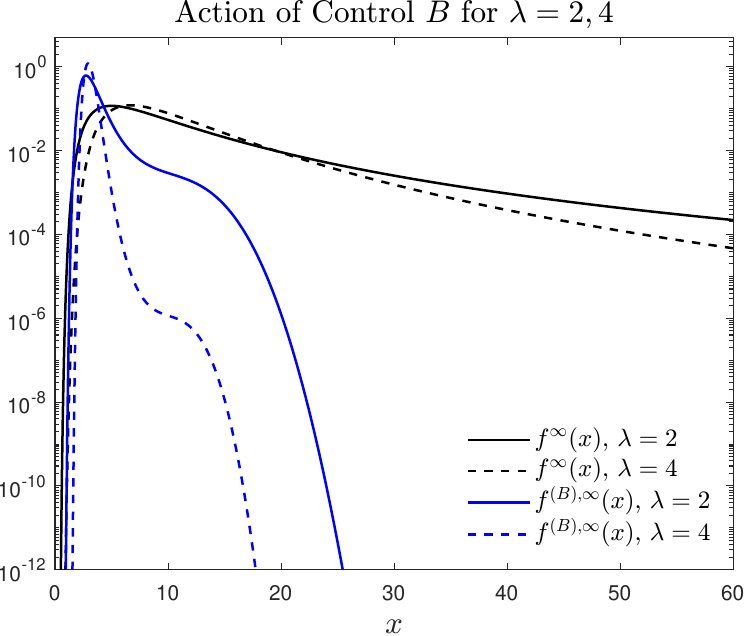}
\caption{Contact distribution at the equilibrium for different coefficients $\lambda=2,4$ and control strategies $A$ (left), and $B$ (right). The other coefficients are: $\delta=-1$, $\sigma^2=0.2$, $m=10$, $x_T=3$, and $\nu=1$.}
\label{fig:test2-3}
\end{figure}

Finally, in Figure \ref{fig:test2-3} we fix the penalization coefficient $\nu=1$ and the vary $\alpha$ with $\sigma^2$ fixed in a way that $\lambda=2,4$. As the value of $\lambda$ decreases, the tails of the distribution increasingly exhibit their power-law behaviour. Consequently, for a fixed penalization coefficient $\nu$, the control reduces them less.

\subsection{Test 3: Consistency in the macroscopic limit without the control}
\begin{figure}[t!]
\centering
\includegraphics[width=0.475\textwidth]{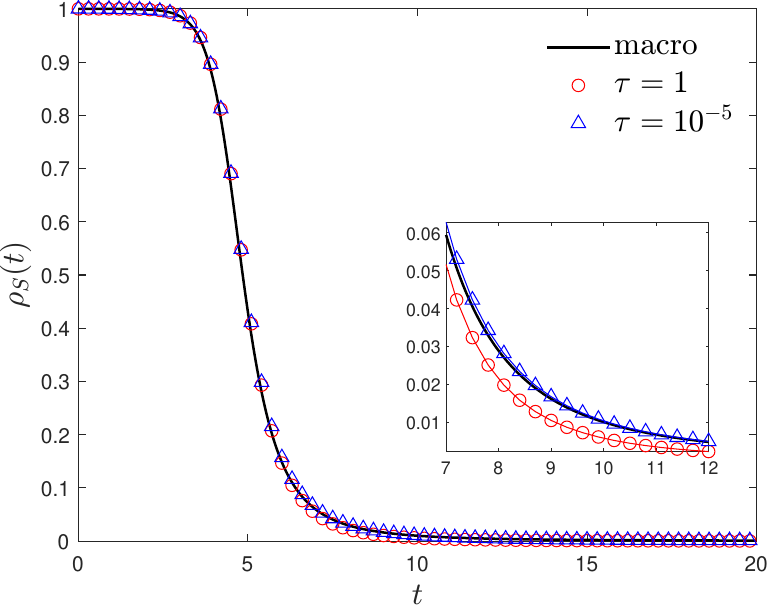}
\includegraphics[width=0.475\textwidth]{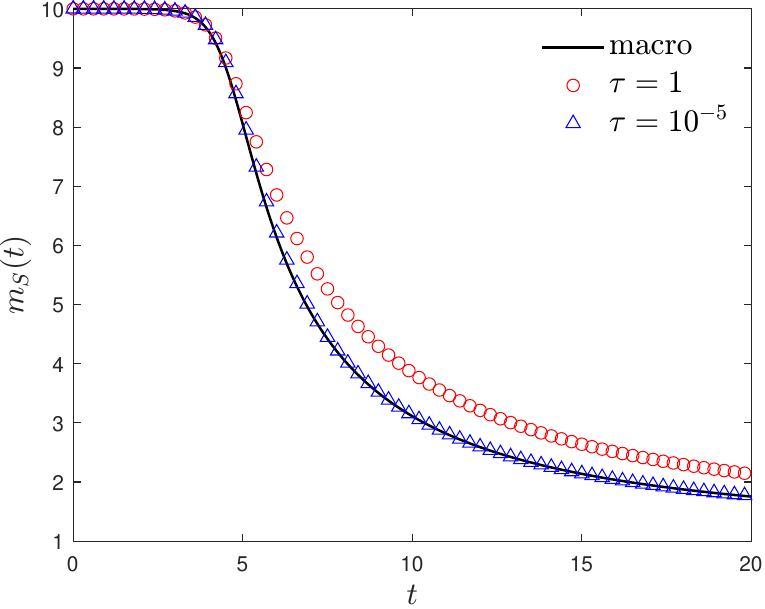} \\
\vspace{5pt}
\includegraphics[width=0.475\textwidth]{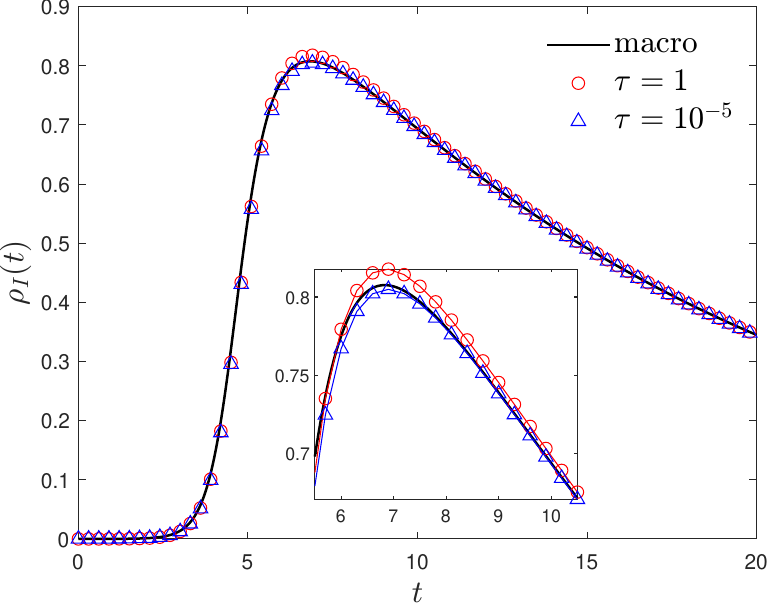}
\includegraphics[width=0.475\textwidth]{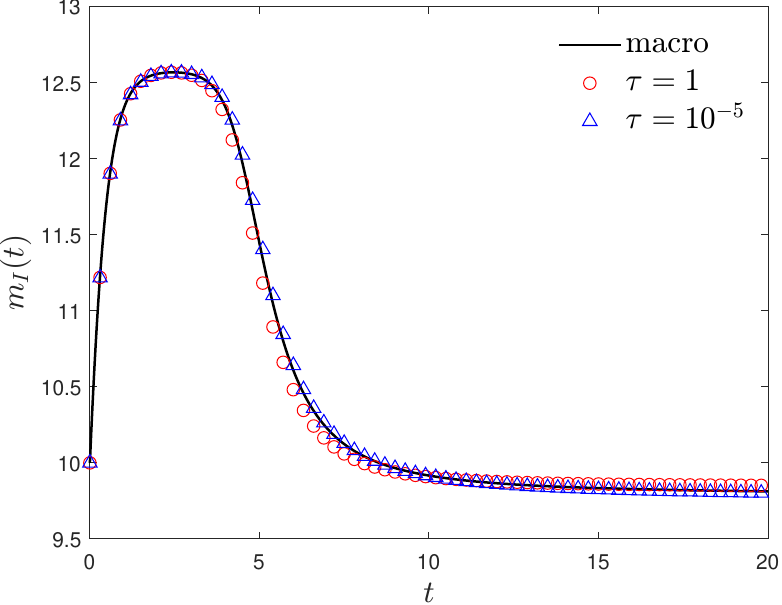}\\
\vspace{5pt}
\includegraphics[width=0.475\textwidth]{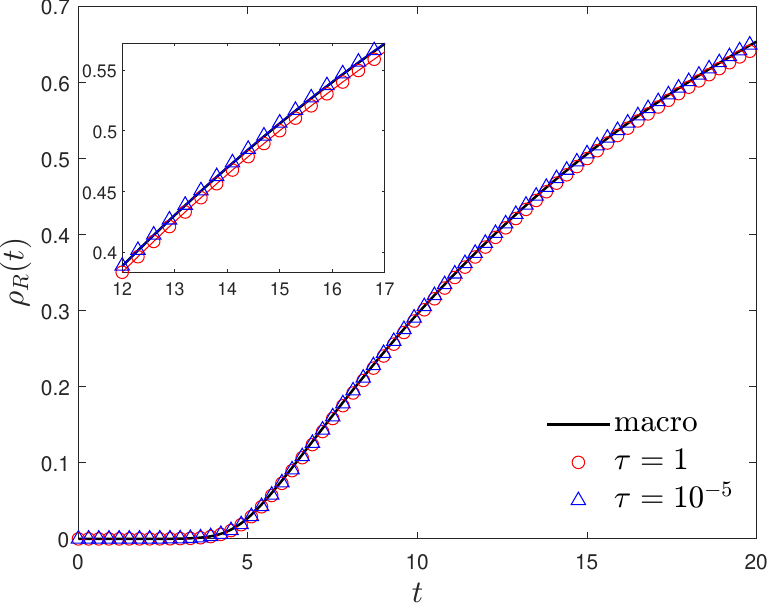}
\includegraphics[width=0.475\textwidth]{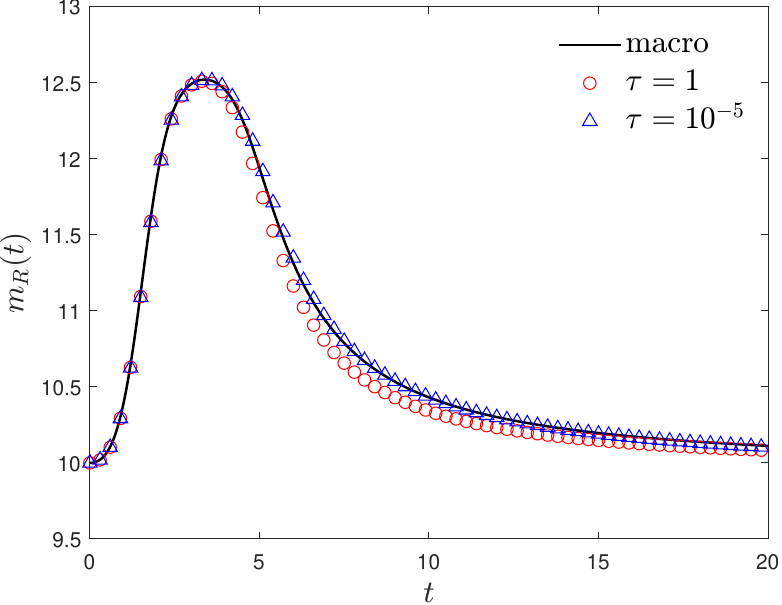}
\caption{Time evolution of the mass fraction (left) and mean (right) of system \eqref{eq:kinetic} with $L=2$, $\tau=1$ (red circles) and $\tau=10^{-5}$ (black triangles), together with the solution of \eqref{eq:mass_L2}-\eqref{eq:mom_L2} (solid black lines). The epidemiological parameters are $\beta_1=2\cdot10^{-2}$, $\beta_2=2\cdot10^{-6}$, $\gamma=1/14$, the other parameters are $\Delta x=0.02$, $\Delta t=0.01$, $\alpha=1$, $\sigma^2=0.2$, in a way that $\lambda=5$, and $\delta=-1$. }
 \label{fig:test3}
\end{figure}

In this test, we investigate numerically the consistency of the macroscopic closure of the kinetic system \eqref{eq:kinetic} with $L=2$ and $\tau\to0$ with the system of equations for the mass and first order moment \eqref{eq:mass_L2}-\eqref{eq:mom_L2}. We choose the $x$-domain $[0,500]$ discretized with step size $\Delta x = 0.02$, while the time domain is $[0,20]$ with time step $\Delta t = 0.01$. The epidemiological parameters are $\beta_1=2\cdot10^{-2}$, $\beta_2=2\cdot10^{-6}$, $\gamma=1/14$, the other parameters are $\alpha=1$, $\sigma^2=0.2$, in a way that $\lambda=5$, and $\delta=-1$. The initial distribution in every compartment $J$ is
\[
f_J(x,0)= \rho_J(0) \frac{\lambda^\lambda}{(m^0_J)^\lambda\Gamma(\lambda)} x^{\lambda-1} \exp\left\{-\frac{\lambda x}{m^0_J}\right\}
\]
with $m^0_J=10$ for every $J$, and initial masses $\rho_I(0)=\rho_R(0)=10^{-5}$ and $\rho_S(0)=1-\rho_I(0)-\rho_R(0)$.

In Figure \ref{fig:test3} we observe that in the limit $\tau \to 0^+$ we obtain a numerical evidence of the consistency between the derived macroscopic system and the kinetic model \eqref{eq:kinetic} with equilibrium closure approximation \eqref{eq:mass_L2}-\eqref{eq:mom_L2}, as expected. 

\begin{figure}
\centering
\includegraphics[width=0.475\textwidth]{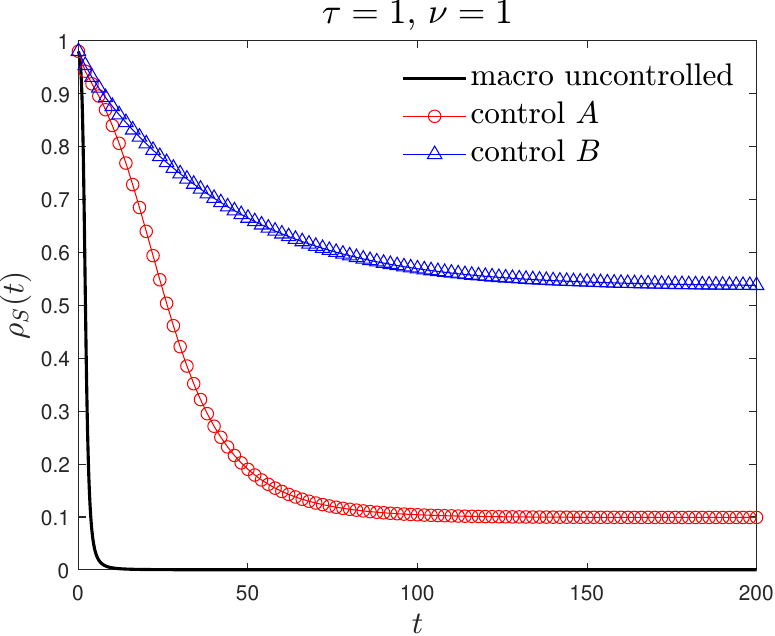}
\includegraphics[width=0.475\textwidth]{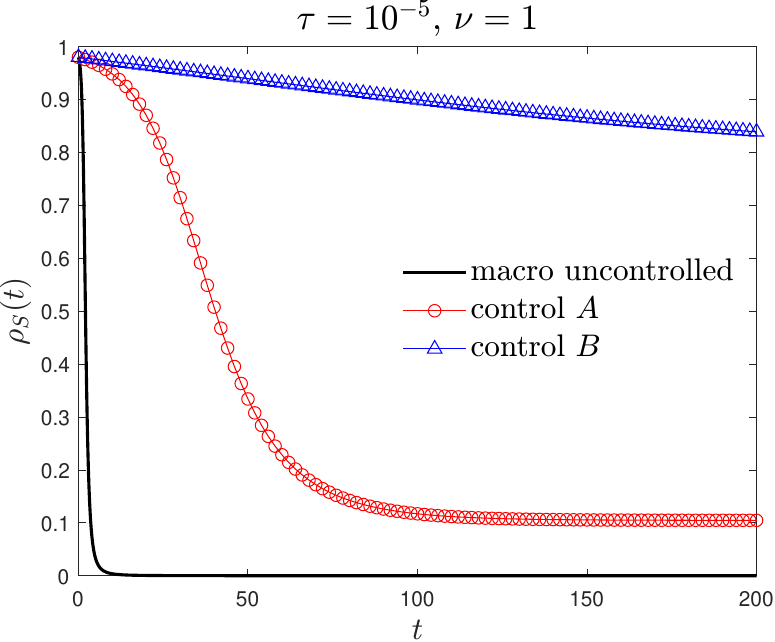} \\
\vspace{5pt}
\includegraphics[width=0.475\textwidth]{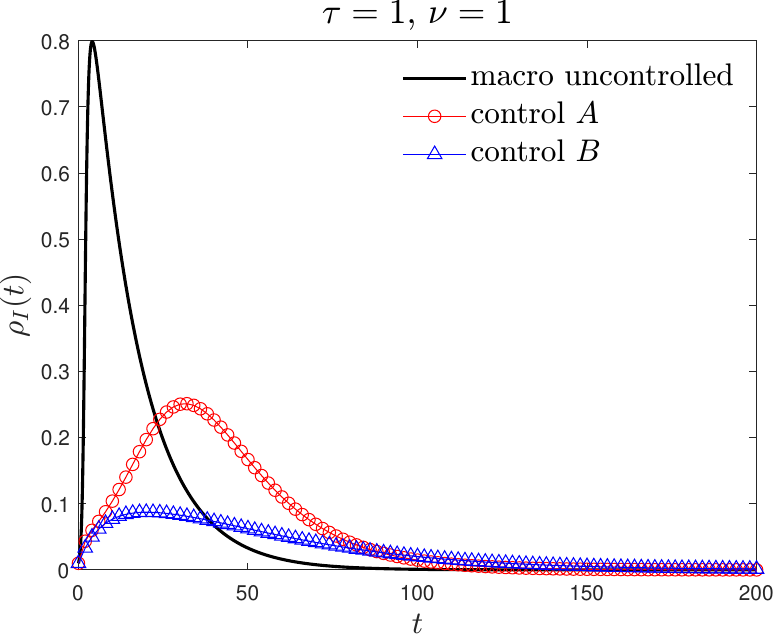}
\includegraphics[width=0.475\textwidth]{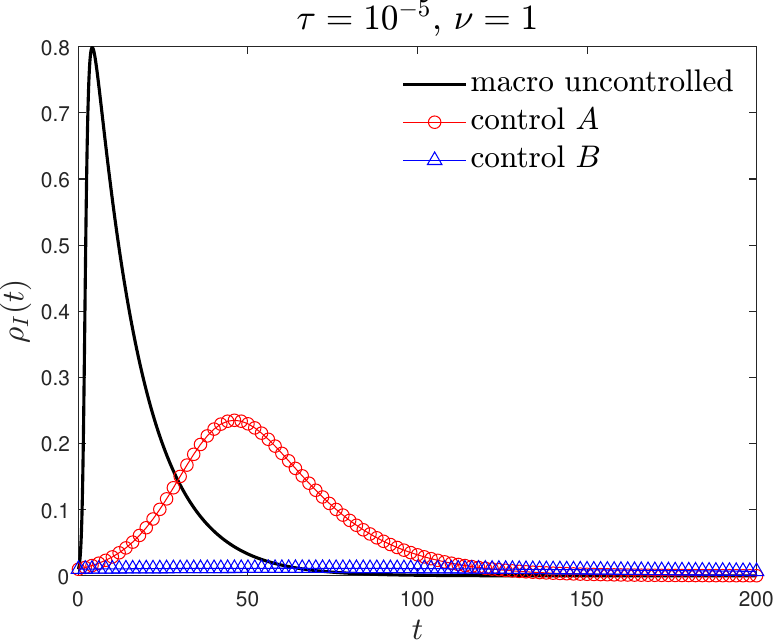}\\
\vspace{5pt}
\includegraphics[width=0.475\textwidth]{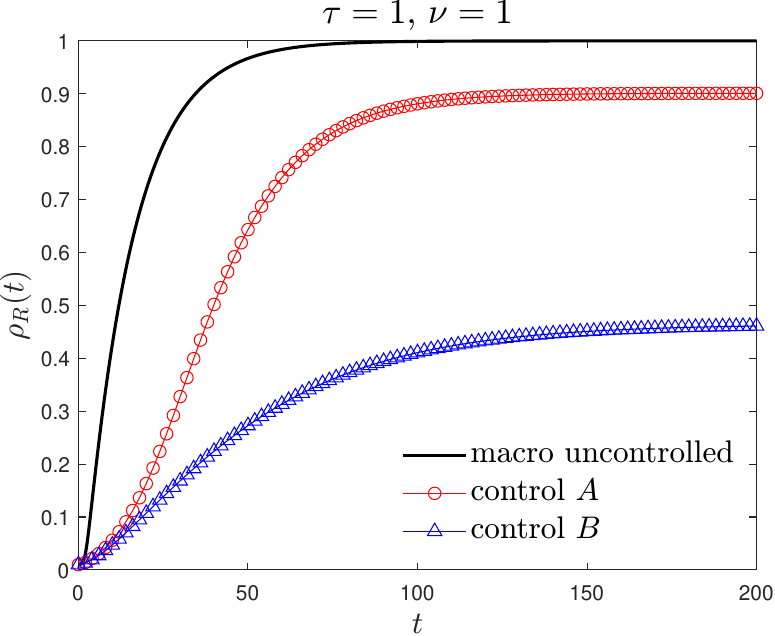}
\includegraphics[width=0.475\textwidth]{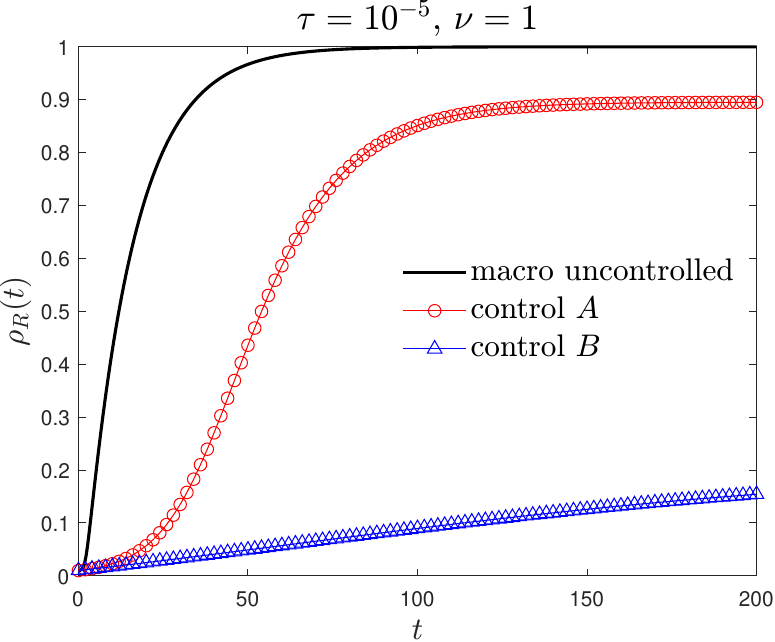}
\caption{Time evolution of the mass fractions $\rho_J(t)$ in every compartment $J=\{S,I,R\}$. Left column: $\tau=1$ and $\nu=1$; right column: $\tau=10^{-5}$ and $\nu=1$. The black lines represent the uncontrolled scenario obtained by solving the system \eqref{eq:mass_L2}-\eqref{eq:mom_L2}. The red circled lines is the control $A$ obtained by solving the system \eqref{eq:kinetic_control} with $H=A$, the black triangulated lines the control $B$, obtained by solving the same system with $H=B$.  The epidemiological parameters are $\beta_1=2\cdot10^{-2}$, $\beta_2=2\cdot10^{-6}$, and $\gamma=1/14$. Initial masses are $\rho_I(0)=\rho_R(0)=10^{-2}$ and $\rho_S(0)=1-2\cdot10^{-2}$. The other parameters are
$\Delta t=0.01$, $\Delta x=0.02$ in the domain $[0,500]$, $\delta=-1$, $\alpha=1$ and $\sigma^2=0.2$.}
\label{fig:test4-1}
\end{figure}

\subsection{Test 4: The controlled system}
In this test, we compare the effect of the different control strategies by numerically solve the system \eqref{eq:kinetic_control} with $H=A,B$. The adopted scheme is described at the beginning of Section \ref{sect:4} with $\Delta t=0.01$ and $\Delta x=0.02$ in the domain $[0,500]$. The epidemiological parameters are $\beta_1=2\cdot10^{-2}$, $\beta_2=2\cdot10^{-6}$, and $\gamma_I=1/14$. We choose in all tests $\delta=-1$, $\alpha=1$ and $\sigma^2=0.2$, so that $\lambda=5$. The initial conditions are
\[
f_J(x,0)= \rho_J(0) \frac{\lambda^\lambda}{(m^0_J)^\lambda\Gamma(\lambda)} x^{\lambda-1} \exp\left\{-\frac{\lambda x}{m^0_J}\right\}
\]
with $m^0_J=10$ for every compartment $J$, and initial masses $\rho_I(0)=\rho_R(0)=10^{-2}$ and $\rho_S(0)=1-\rho_I(0)-\rho_R(0)$.

In Figure \ref{fig:test4-1}, we solve the kinetic system \eqref{eq:kinetic_control} for several values of $\tau >0$ and we compare the results with the unconstrained system with $L = 2$. In this test, we fix the penalization coefficient $\nu=1$ and we compare the effect of the control strategies $(A)$ and $(B)$ for different scale parameter $\tau=1,10^{-5}$. We may observe how the control $(B)$, while shaping the contact distribution in a slim-tailed distribution, is able to impact on the infection dynamics. On the other hand, the control $(A)$  as it does not modify the structure of the contact distribution, which retains its overpopulated tail. 


\section*{Conclusion}
In this work, we focused on the macroscopic effects of kinetic control strategies applied to epidemic dynamics. We discussed the derivation of \textcolor{black}{general macroscopic epidemic models, in which the epidemic transmission dynamics is influenced by a general number of moments of the contact distribution. These models describe the evolution of a controlled multi-agent system, obtained through an equilibrium closure approach.} In more detail, we considered two control strategies previously proposed in the literature for similar models. Our results demonstrated how controlling the interaction strengths can transform a contact distribution with an overpopulated tail into a slim-tailed distribution, thereby reducing the heterogeneity in interactions. Additionally, we discussed the role of power-law tails in epidemic dynamics within the framework of the derived macroscopic models. Numerical results highlighted the effectiveness of the proposed approach in mitigating the spread of infection and demonstrated its potential for practical applications in epidemic control.

\section*{Acknowledgments}
The work has been written within the activities of GNFM group of INdAM (National Institute of High Mathematics). 
A.M. was supported by the Advanced Grant Nonlocal-CPD (Nonlocal PDEs for Complex Particle Dynamics: Phase Transitions, Patterns and Synchronization) of the European Research Council Executive Agency (ERC) under the European Union’s Horizon 2020 research and innovation programme (grant agreement No.~883363). M.Z. acknowledges partial support by the MUR-PRIN 2020 project (No. 2020JLWP23) “Integrated Mathematical Approaches to Socio–Epidemiological Dynamics” and by European Union - NextGeneration EU.


\end{document}